\documentclass[11pt]{article}
\usepackage{color,enumerate}
\usepackage{graphics}
\usepackage{graphics}
\usepackage{tikz}
\usepackage[all]{xy}
\usepackage{amssymb}
\definecolor{green}{rgb}{0,.6,0}
\definecolor{purple}{rgb}{0.5,0,0.5}

\def\n{\mathfrak n}

\newfont{\bbb}{msbm10 scaled\magstephalf}
\newfont{\sbbb}{msbm7 scaled\magstephalf}
\def\C{\mbox{\bbb{C}}}
\def\P{\mbox{\bbb{P}}}
\def\R{\mbox{\bbb{R}}}
\def\Z{\mbox{\bbb{Z}}}
\def\SZ{\mbox{\sbbb{Z}}}
\def\SC{\mbox{\sbbb{C}}}
\def\cd{\C^d}
\def\rd{\R^d}
\def\rn{\R^n}
\def\zd{\Z^d}
\def\td{T^d}
\def\rddu{(\rd)^*}
\def\rndu{(\rn)^*}
\def\et1{e^{2\pi i\theta_1}}
\def\etd{e^{2\pi i\theta_d}}

\def\vz{\underline{z}}
\def\vw{\underline{w}}
\def\zjs{|z_j|^2}

\def\D{\Delta}
\def\xd{X_1,\ldots,X_d}
\def\ld{\lambda_1,\ldots,\lambda_d}
\def\pic{\pi_{\SC}}
\newtheorem{thm}{Theorem}[section]
\newtheorem{defn}[thm]{Definition}
\newtheorem{remark}[thm]{Remark}

\def\squareforqed{\hbox{\rlap{$\sqcap$}$\sqcup$}}
\def\qed{\ifmmode\else\unskip\quad\fi\squareforqed}
\def\smartqed{\def\qed{\ifmmode\squareforqed\else{\unskip\nobreak\hfil
\penalty50\hskip1em\null\nobreak\hfil\squareforqed
\parfillskip=0pt\finalhyphendemerits=0\endgraf}\fi}}

\newcounter{sect}

\title{\sc Toric Geometry of the Regular Convex Polyhedra}
\author{\sc Fiammetta Battaglia and Elisa Prato
\vspace{0.3cm}\\
\sc Dipartimento di Matematica e Informatica "U. Dini"\\
\sc Universit\`a di Firenze \\
\sc Viale Morgagni,  67/A\\
\sc 50134 Firenze, ITALY \\
{\tt fiammetta.battaglia@unifi.it}
\sc and
{\tt elisa.prato@unifi.it}}
\date{}
\begin{document}
\maketitle
\begin{abstract}
In this article, we describe symplectic and complex toric spaces associated to the five regular convex polyhedra.
The regular tetrahedron and the cube are rational and simple, the regular octahedron is not simple, the regular dodecahedron is not rational and the regular icosahedron is neither simple nor rational. We remark that the last two cases cannot be treated via standard toric geometry.
\end{abstract}

\section*{Introduction}
Among the five regular convex polyhedra, the regular tetrahedron and the cube are examples of simple rational convex polyopes. To these, the standard smooth toric geometry applies, both in the symplectic and complex category. From the symplectic viewpoint, in fact, the regular tetrahedron and the cube satisfy the hypotheses of the Delzant Theorem \cite{d}, and it is easily seen that they correspond, respectively, to $S^7/S^1$ and $S^2\times S^2 \times S^2$. From the complex viewpoint, on the other hand, the toric variety associated to the regular tetrahedron is $\C\P^3$ while $\C\P^1\times\C\P^1\times\C\P^1$ corresponds to the cube (see for example \cite[Section 1.5]{fulton}); these can also be obtained as quotients by a complex version of the Delzant procedure described by Audin in \cite[Chapter VI]{au}.

The regular octahedron is still rational, but it is no longer simple. The toric variety associated to the regular octahedron is however well known and is described, for example, by Fulton in \cite[Section 1.5]{fulton}; it can also be obtained as a complex quotient by applying the Cox construction \cite[Theorem~2.1]{cox}.

The regular dodecahedron, on the other hand, is simple but it is the first of the five regular convex polyhedra that is not rational. It is shown by Prato in \cite{pdoc} that, by applying her extension of the Delzant procedure to the case of general simple convex polytopes \cite{p}, one can associate to the regular dodecahedron a {\em symplectic toric quasifold}. Quasifolds are a generalization of manifolds and orbifolds: they are not necessarily Hausdorff and they are locally modeled by the quotient of a manifold modulo the action of a discrete group. 

In this article, we recall all of the above and we complete the picture, first of all, by associating to the regular dodecahedron a {\em complex toric quasifold}. We do so by applying a generalization, given by the authors in \cite{cx}, of the 
procedure described by Audin in \cite[Chapter VI]{au} to the case of general simple convex polytopes. 
As in the smooth case, the symplectic and complex quotients can be identified \cite[Theorem~3.2]{cx}, 
endowing the corresponding toric quasifold with a K\"ahler structure.

We go on to address the case of the regular icosahedron. From the toric viewpoint, this is certainly the most complicated of the five regular convex polyhedra, since it is neither simple nor rational. However, we can apply generalizations by Battaglia and Prato \cite{stratificato-announ} and Battaglia \cite{stratificato-re, stratificato-cx} of toric quasifolds \cite{p,cx} and of the Cox construction \cite{cox} to arbitrary convex polytopes; this allows us to associate to the regular icosahedron, both in the symplectic and complex category, a space that is stratified by quasifolds. As for all $3$--dimensional polytopes, here there are only 
zero--dimensional singular strata and an open dense regular stratum. Moreover, by \cite[Theorem~3.3]{stratificato-cx}, the symplectic and complex quotients can be identified, endowing the regular stratum with the structure of a K\"ahler quasifold.

Notice, finally, that we are still missing a symplectic toric space corresponding to the regular octahedron; this too can be found by applying Battaglia's work on arbitrary convex polytopes. What we get here is a space that is stratified by symplectic manifolds (see \cite[Remark~6.6]{stratificato-re}); moreover, by \cite[Theorem~3.3]{stratificato-cx}, this symplectic quotient can be identified with the complex quotient, and the regular stratum is K\"ahler.

The article is structured as follows: in Section~\ref{polytopes} we recall a few necessary facts on convex polytopes; in Section \ref{simple} we recall 
from \cite{p,cx} how to construct symplectic and complex toric quasifolds from simple convex polytopes; in Section~\ref{nonsimple} we recall from \cite{stratificato-re, stratificato-cx} the construction of the symplectic and complex toric spaces corresponding to arbitrary convex polytopes; 
finally, in Sections~\ref{smooth}, \ref{octahedron}, \ref{dodecahedron}, and \ref{icosahedron} we describe the symplectic and complex toric 
spaces corresponding to the five regular convex polyhedra.

\section{Facts on convex polytopes}\label{polytopes}
Consider a dimension $n$ convex polytope $\D\subset\rndu$.
\begin{defn}[Simple polytope]{\rm $\D$ is said to be {\em simple} if each of its vertices is contained in exactly $n$ facets.}\end{defn} 
Assume now that $\Delta$ has $d$ facets. Then there exist elements $\xd$ in $\rn$ and real numbers $\ld$ such that
\begin{equation}\label{polydecomp}
\D=\bigcap_{j=1}^d\{\;\mu\in\rndu\;|\;\langle\mu,X_j\rangle\geq\lambda_j\;\}.
\end{equation}
Let us consider the open faces of $\Delta$. They can be described
as follows. For each such face $F$ there exists a, possibly empty,
subset $I_F\subset\{1,\ldots,d\}$ such that
\begin{equation}\label{facce}
F=\{\,\mu\in\Delta\;|\;\langle\mu,X_j\rangle=\lambda_j\;
\hbox{ if and only if}\; j\in I_F\,\}.
\end{equation}
A partial order on the set of all open faces of $\D$ is defined by
setting $F\leq F'$ (we say $F$ contained in $F'$) if
$F\subseteq\overline{F'}$. Notice that $F\leq F'$ if, and only if, $I_{F'}\subset I_F$.
The polytope $\Delta$ is the disjoint
union of its open faces. Let $r_F=\hbox{card}(I_F)$; we have the
following definitions:
\begin{defn}{\rm A $p$--dimensional open face $F$ of the polytope is
said to be {\em singular} if $r_F>n-p$.}
\end{defn}
\begin{defn}{\rm A $p$--dimensional open face $F$ of the polytope is
said to be {\em regular} if $r_F=n-p$.}
\end{defn}
\begin{remark}{\rm Let $F$ be a $p$--dimensional singular face in $(\R^n)^*$, 
then $p<n-2$.
Therefore any polytope in $(\R^2)^*$ is simple and the singular faces of a 
nonsimple polytope in $(\R^3)^*$ are vertices. 
}\end{remark}
We refer the reader to Ziegler's book \cite{ziegler} for additional basic facts on convex polytopes.
We now go on to recall what is meant by quasilattice, and quasirational polytope.
\begin{defn}[Quasilattice]{\rm
A {\em quasilattice} in $\rn$ is the $\Z$--span of a set of $\R$--spanning vectors, $Y_1,\ldots,Y_q$, of
$\rn$.}
\end{defn}
Notice that $\hbox{Span}_{\SZ}\{Y_1,\dots,Y_q\}$ is a lattice if, and
only if, it is generated by a basis of $\rn$.
\begin{defn}[Quasirational polytope]{\rm Let $Q$ be a quasilattice in
$\rn$. A convex polytope $\D\subset\rndu$ is said to be {\em
quasirational} with respect to the quasilattice $Q$ if the vectors $\xd$ in
(\ref{polydecomp}) can be chosen in $Q$.}
\end{defn}
Remark that each polytope in $\rndu$ is quasirational with
respect to the quasilattice $Q$ that is generated by the elements $\xd$ in (\ref{polydecomp}). 
We note that if $\xd$ can be chosen inside a lattice, then the polytope is rational.

\section{The simple case}\label{simple}
Let $\Delta\subset \rndu$ be an $n$--dimensional simple convex polytope. 
We are now ready to recall from \cite{p} and \cite{cx} the construction of the symplectic and complex toric quasifolds associated to $\Delta$.
For the definition and main properties of symplectic and complex quasifolds  we refer the reader to \cite{p,kite} and \cite{cx}.
For the purposes of this article, we will restrict our attention to the special case $n=3$.
We begin by remarking that both constructions rely on the notion of quasitorus, which we recall.
\begin{defn}[Quasitorus]{\rm Let $Q$ be a quasilattice in $\R^3$. 
We call {\em quasitorus} of dimension $3$ the group and quasifold $D^3=\R^3/Q$.}
\end{defn}
Notice that, if the quasilattice is a lattice, we obtain the classical notion of torus. 
The quasilattice $Q$ also acts naturally on $\C^3$:
$$
\begin{array}{lcccc}
Q&\times&\C^3 & \longrightarrow & \C^3\\ (A&,&X+iY) & \longmapsto & (X+A)+iY.
\end{array}
$$
Therefore, in the complex category we have
\begin{defn}[Complex quasitorus]
{\rm Let $Q$ be a quasilattice in $\R^3$. We call {\em complex quasitorus} of dimension $3$ the group and 
complex quasifold $D^3_{\SC}=\C^3/Q$.}
\end{defn}
In analogy with the smooth case, we will say that $D^3_{\SC}$ is the {\em complexification} of $D^3$.
Assume now that our polytope $\Delta$ is quasirational with respect to a quasilattice $Q$ and write 
\begin{equation}\label{decomp}
\D=\bigcap_{j=1}^d\{\;\mu\in(\R^3)^*\;|\;\langle\mu,X_j\rangle\geq\lambda_j\;\}
\end{equation}
for some elements $\xd\in Q$ and some real numbers $\ld$; again, $d$ here is the number of facets of $\Delta$.
Let $\{e_1,\ldots,e_d\}$ denote the standard basis
of $\R^d$ and $\C^d$. Consider the surjective linear mapping
$$
\begin{array}{cccc}\label{pi}
\pi \,\colon\,& \R^d & \longrightarrow & \R^3\\
    &   e_j& \longmapsto & X_j,
\end{array}
$$ and its complexification $$
\begin{array}{cccc}
\pic \,\colon\,& \C^d & \longrightarrow & \C^3\\
    &   e_j& \longmapsto & X_j.
\end{array}
$$ 
Consider the quasitori $D^3=\R^3/Q$ and $D^3_{\SC}=\C^3/Q$. The mappings $\pi$ and $\pic$ each induce group epimorphisms
$$\Pi\,\colon\, T^d=\R^d/\Z^d\longrightarrow D^3$$ 
and 
$$\Pi_{\SC} \,\colon\, T^d_{\SC}=\C^d/\Z^d\longrightarrow D^3_{\SC}.$$
We define $N$ to be the kernel of the mapping $\Pi$ and $N_{\SC}$ to be the kernel of the
mapping $\Pi_{\SC}$. Notice that neither $N$ nor $N_{\SC}$ are honest tori unless $Q$ is a honest
lattice. The Lie algebras of $N$ and $N_{\SC}$ are, respectively, $\n=\ker \pi$ and $\n_{\SC}=\ker\pic$.
The mappings $\Pi$ and $\Pi_{\SC}$ induce isomorphisms
$$
T^d/N\longrightarrow D^3
$$
and
$$
T^d_{\SC}/N_{\SC}\longrightarrow D^3_{\SC}.
$$
Let us begin with the symplectic construction.
Consider the space $\cd$, endowed with the symplectic form $\omega_0=\frac{1}{2\pi
i}\sum_{j=1}^d dz_j\wedge d\bar{z}_j$ and the action of the torus $\td=\rd/\zd$:
$$
\begin{array}{cccccl}
& \td&\times&\cd&\longrightarrow& \cd\\
&((\et1,\ldots,\etd)&,&\vz)&\longmapsto&(\et1 z_1,\ldots, \etd z_d).
\end{array}
$$
This action is effective and Hamiltonian, with moment mapping given by
$$
\begin{array}{cccl}
J\,\colon&\cd&\longrightarrow &\rddu\\
&\vz&\longmapsto & \sum_{j=1}^d \zjs
e_j^*+\lambda,\quad\lambda\in\rddu \;\mbox{constant}.
\end{array}
$$
Choose now $\lambda=\sum_{j=1}^d {\lambda_j} e_j^*$, with $\ld$ as in (\ref{decomp}).
Denote by $i$ the Lie algebra inclusion $\n\rightarrow\rd$ and
notice that $\Psi=i^*\circ J$ is a moment mapping for the induced
action of $N$ on $\cd$.  Consider now the orbit space $M_{\Delta}=\Psi^{-1}(0)/N$.
Then we have, from \cite[Theorem 3.3]{p}:
\begin{thm}[Generalized Delzant construction]\label{thmp1}
Let $Q$ be a quasilattice in $\R^3$ and let $\D\subset(\R^3)^*$ be a $3$--dimensional simple convex polytope
that is quasirational with respect to $Q$. Assume that $d$ is the number of facets of $\D$ and
consider vectors $X_1,\ldots,X_d$ in $Q$ that satisfy (\ref{decomp}). For each $(\D, Q,\{\xd\})$,
the orbit space $M_{\Delta}$ is a compact, connected $6$--dimensional symplectic quasifold
endowed with an effective Hamiltonian action of the quasitorus $D^3=\R^3/Q$ such that, 
if $\Phi\,\colon M\rightarrow(\R^3)^*$ is the corresponding moment mapping, then $\Phi(M_{\Delta})=\Delta$. 
\end{thm}
We say that the quasifold $M_{\Delta}$ with the effective Hamiltonian action of $D^3$ is the {\em symplectic toric quasifold} 
associated to $(\D, Q, \{\xd\})$.

Let us now pass to the complex construction. Following the notation of the previous section, consider, for any open face $F$ of
$\Delta$, the $T^d_{\SC}$--orbit
$$\C^d_F=\{\,\underline{z}\in\cd\;|\;z_j=0\;\;\hbox{iff}\;\; j\in I_F\,\}.$$ 
Consider the open subset of $\C^d$ given by
\begin{equation}\label{apertopulito}
\C^d_{\D}=\bigcup_{F\in\D} \{\,\underline{z}\in\cd\;|\;z_j\neq0\;\;\hbox{if}\;\; j\notin I_{F}\,\}.
\end{equation}
Notice that \begin{equation}\label{aperto}
\C^d_{\D}=\bigcup_{\mu\in\D} \{\,\underline{z}\in\cd\;|\;z_j\neq0\;\;\hbox{if}\;\; j\notin I_{\mu}\,\},
\end{equation} 
where $\mu$ ranges over all the vertices of the polytope $\Delta$.
Moreover, since the polytope is simple, we have that $$\C^d_{\Delta}=\bigcup_F\C^d_F.$$
In fact, in this case,
$\{\,\underline{z}\in\cd\;|\;z_j\neq0\;\;\hbox{if}\;\; j\notin I_{F}\,\}=\cup_{F'\geq F}\C^d_{F'}$.
The group $N_{\SC}$ acts on the space $\C^d_{\Delta}$.
Consider the space of orbits $X_{\Delta}=\C^d_{\Delta}/N_{\SC}$. We then have, from \cite[Theorem 2.2]{cx}:
\begin{thm}\label{poltocx} 
Let $Q$ be a quasilattice in $\R^3$ and let $\D\subset(\R^3)^*$ be a $3$--dimensional simple convex polytope
that is quasirational with respect to $Q$. Assume that $d$ is the number of facets of $\D$ and
consider vectors $X_1,\ldots,X_d$ in $Q$ that satisfy (\ref{decomp}). For each $(\D, Q,\{\xd\})$,
the corresponding quotient $X_{\Delta}$ is a complex quasifold of dimension $3$, 
endowed with a holomorphic action of the complex quasitorus $D^3_{\SC}=\C^3/Q$ having a dense open orbit.
\end{thm}
We say that the quasifold $X_{\Delta}$ with the holomorphic action of $D^3_{\SC}$ is the {\em complex toric quasifold} 
associated to $(\D, Q, \{\xd\})$.

Finally, we conclude this section by recalling that the natural embedding
\begin{equation}
\Psi^{-1}(0)\hookrightarrow\C^d_{\Delta},
\end{equation}
induces a mapping
$$
\chi\,\colon\, M_{\Delta} \longrightarrow X_{\Delta}
$$
that sends each $N$--orbit to the corresponding $N_{\SC}$--orbit. This mapping is equivariant
with respect to the actions of the quasitori $D^3$ and $D^3_{\SC}$. 
Then, under the same assumptions of Theorems \ref{thmp1} and \ref{poltocx}, we have, from \cite[Theorem 3.2]{cx}:
\begin{thm}
\label{teoremadellachi} 
The mapping $\chi\,\colon\,M_{\Delta} \longrightarrow X_{\Delta}$ is an equivariant
diffeomorphism of quasifolds. Moreover, the induced symplectic form on the complex
quasifold $X_{\Delta}$ is K\"ahler.
\end{thm}

For the smooth case see Audin \cite[Proposition~3.1.1]{au} but also Guillemin \cite[Appendix 1, Theorem~1.4]{vg}.
\section{The nonsimple case}\label{nonsimple}
Consider now a nonsimple convex polytope $\Delta\subset(\R^3)^*$ and assume that $\Delta$ is quasirational with respect to a quasilattice $Q$.
The idea here is to repeat the constructions of the previous section. If we do so, we again find the groups $N$ and $N_{\SC}$ and the quasitori $D^3$ and $D^3_{\SC}$, isomorphic to  $T^d/N$ and $T^d_{\SC}/N_{\SC}$ respectively.  However, the symplectic construction produces spaces that are {\em stratified by symplectic quasifolds}, while the complex construction yields spaces that are {\em stratified by complex quasifolds}.
For the exact definitions of these notions we refer the reader to \cite[Section~2]{stratificato-re} and \cite[Definition~1.5]{stratificato-cx}. We remark that in the smooth case these definitions yield the classical definition of Goresky and MacPherson \cite{gmp}.
Many of the important features of these stratified structures will be clarified when addressing the relevant examples (see Sections \ref{octahedron} and \ref{icosahedron}).

Let us consider the symplectic case first. The main difference with respect to the case of simple polytopes is that here there are points in the level set
$\Psi^{-1}(0)$ that have isotropy groups of positive dimension; therefore $\Psi^{-1}(0)$ is no longer a smooth manifold. 
From Proposition~3.3 and Theorems 5.3, 5.10, 5.11, 6.4 in \cite{stratificato-re} we have
\begin{thm}[Generalized Delzant construction: nonsimple case]\label{nonsemplice-re}
Let $Q$ be a quasilattice in $\R^3$ and let $\D\subset(\R^3)^*$ be a $3$--dimensional convex polytope
that is quasirational with respect to $Q$. Assume that $d$ is the number of facets of $\D$ and
consider vectors $X_1,\ldots,X_d$ in $Q$ that satisfy (\ref{decomp}). For each $(\D, Q,\{\xd\})$,
the quotient $M_{\Delta}$ is a compact, connected $6$--dimensional space stratified by symplectic quasifolds,
endowed with an effective continuos action of the quasitorus $D^3=\R^3/Q$. Moreover, there exists a continuos mapping 
$\Phi\,\colon M\rightarrow(\R^3)^*$ such that $\Phi(M_{\Delta})=\Delta$. Finally, the restriction of the $D^3$--action 
to each stratum is smooth and Hamiltonian, with moment mapping given by the restriction of $\Phi$. 
\end{thm}
When the polytope is rational, these quotients are examples of the symplectic stratified spaces described by Sjamaar--Lerman in \cite{sl}; 
in particular, the strata are either manifolds or orbifolds \cite[Remark~6.6]{stratificato-re}.
We remark that the nonsimple rational case was addressed also by Burns--Guillemin--Lerman in \cite{bgl,bgl1}; 
they gave an in--depth treatment and a classification theorem in the case of isolated singularities.

Let us now examine the complex case. Here, one still considers the open subset $\C^d_{\Delta}$ as defined in (\ref{apertopulito})
but, while in the simple case the orbits of $\exp(i\n)$ on $\C^d_{\Delta}$ were closed, here there are nonclosed orbits. 
We recall first that $N_{\SC}=NA$,
where $A=\exp(i \n)$; this actually happens also in the simple case.
Then, from \cite[Theorem~2.1]{stratificato-cx}, we have
\begin{thm}[Closed orbits]\label{closedorbits}
Let $\vz\in\C^d_{\D}$. Then the $A$--orbit through $\vz$, $A\cdot\vz$, is 
closed if, and only if, there exists a face $F$ such that $\vz$
is in $\C^d_F$.
Moreover, if $A\cdot\vz$ is nonclosed, then its closure contains one, and only one,
closed $A$--orbit.
\end{thm}
Therefore, in order to define a notion of quotient, one defines the following equivalence relation:
two points $\vz$ and $\vw$ in $\C^d_{\Delta}$ are equivalent
with respect to the action of the group $N_{\SC}$ if, and only if,
$$\left(N(\overline{A\cdot\vz})\right)\cap\left(\overline{A\cdot\vw}\right)\neq\emptyset,
$$
where the closure is meant in $\C^d_{\D}$. The space $X_{\D}$ is then defined to be the quotient $\C^d_{\D}//N_{\SC}$ with respect this
equivalence relation. Notice that, if the polytope is simple, $\C^d_{\D}=\cup_{F\in\D}\C^d_F$, thus, by Theorem~\ref{closedorbits}, 
$A$--orbits through points in $\C^d_{\Delta}$ are always closed and the quotient $X_{\D}$ is just the orbit space endowed with the quotient 
topology. From \cite[Proposition~3.1,Theorem~3.2]{stratificato-cx} we have
\begin{thm}\label{nonsemplice-cx} 
Let $Q$ be a quasilattice in $\R^3$ and let $\D\subset(\R^3)^*$ be a $3$--dimensional convex polytope
that is quasirational with respect to $Q$. Assume that $d$ is the number of facets of $\D$ and
consider vectors $X_1,\ldots,X_d$ in $Q$ that satisfy (\ref{decomp}). For each $(\D, Q,\{\xd\})$,
the corresponding quotient $X_{\Delta}$ is endowed with a stratification by complex quasifolds of dimension $3$.
The complex quasitorus $D^3_{\SC}$ acts continuously on
$X_{\Delta}$, with a dense open orbit. Moreover, the restriction of the $D^3_{\SC}$--action to each stratum
is holomorphic. 
\end{thm}
We remark that, when $Q$ is a lattice and the vectors $X_1,\ldots,X_d$ are primitive in $Q$, the quotient $X_{\Delta}$
coincides with the Cox presentation \cite{cox} of the classical toric variety that corresponds to the fan normal to the polytope $\Delta$.
As for classical toric varieties, there is a one--to--one correspondence between
$p$--dimensional orbits of the quasitorus $D_{\SC}^3$ and $p$--dimensional faces of the polytope.
In particular, the dense open orbit corresponds to the interior of the polytope and
the singular strata correspond to singular faces.

We remark that, like in the simple case, the natural embedding $\Psi^{-1}(0)\hookrightarrow\C^d_{\D}$ induces an identification
between symplectic and complex quotients. From \cite[Theorem~3.3]{stratificato-cx} we have
\begin{thm}\label{non semplice iso} 
The mapping $\chi_{\D}\,\colon\,M_{\Delta}\longrightarrow X_{\Delta}$ is a
homeomorphism which is equivariant with respect to the actions of $D^3$ and $D_{\SC}^3$, respectively.
Moreover, the restriction of $\chi_{\D}$ to each stratum is a diffeomorphism of quasifolds. 
Finally, the induced symplectic form on each stratum
is compatible with its complex structure, so that each stratum is K\"ahler.
\end{thm}
We conclude by pointing out that $M_{\Delta}\simeq X_{\D}$ has two different kinds of 
singularities, namely the stratification and the quasifold structure of the strata. The nonsimplicity of the polytope 
yields the decomposition in strata of the corresponding
topological space, whilst its nonrationality produces the quasifold structure of the
strata and also intervenes in the way the strata are glued to each other. 
This last feature can be observed only in spaces with strata of positive dimension; 
this led to a definition of stratification that naturally extends the usual one \cite[Section~2]{stratificato-re}.

\section{Simple and rational: the regular tetrahedron and the cube}\label{smooth}
The regular tetrahedron and the cube are both simple and rational.
Let us recall the construction of the corresponding symplectic and complex toric manifolds.
We follow the notation of Section~\ref{simple}, which also applies to the smooth case.

\begin{figure}[h]
\begin{minipage}[b]{0.45\linewidth}
\begin{center}
\includegraphics[width=50mm]{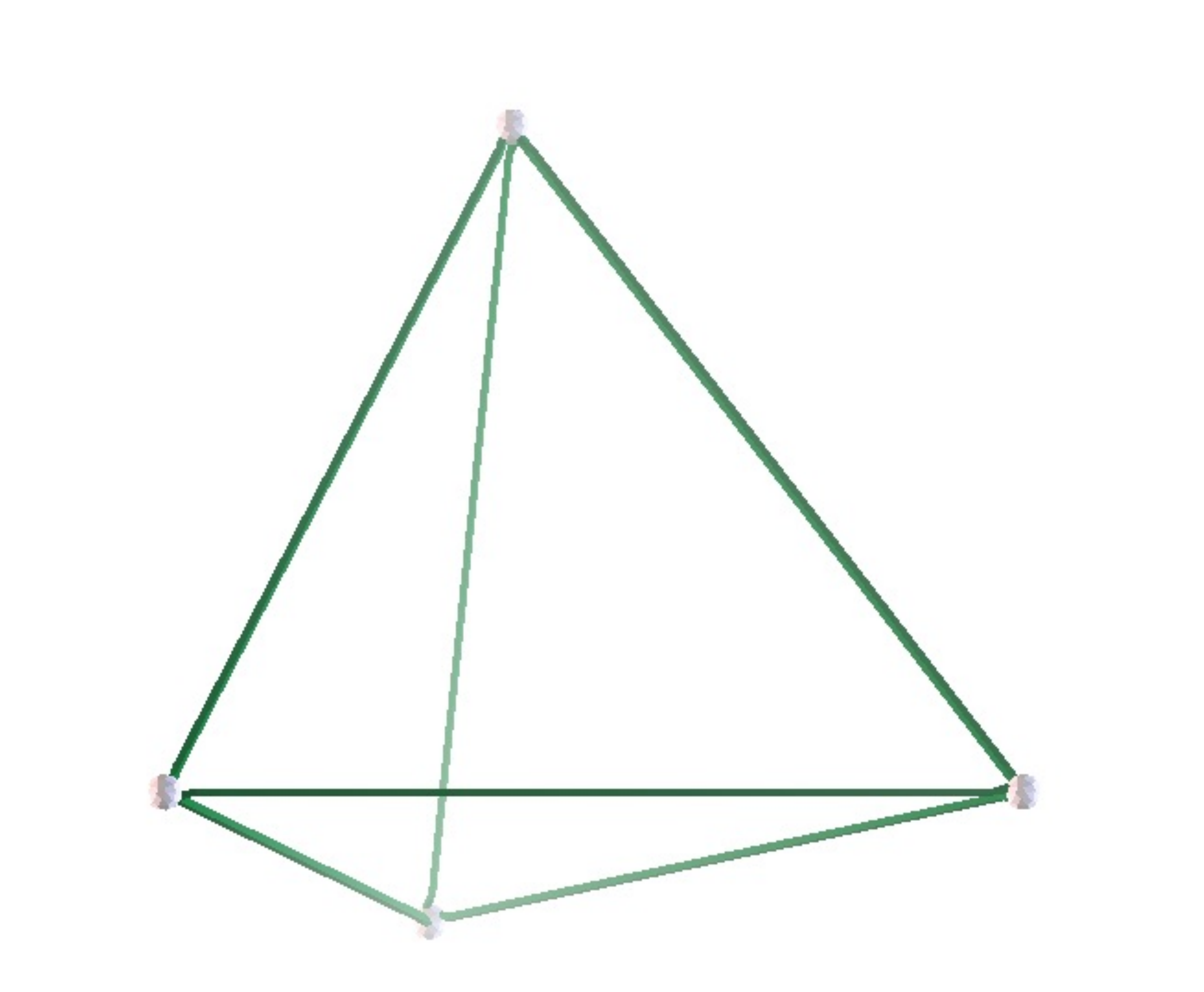}
\end{center}
\qquad
\caption{The regular tetrahedron}\label{tetraedro}
\end{minipage}
\begin{minipage}[b]{0.45\linewidth}
\begin{center}
\includegraphics[width=50mm]{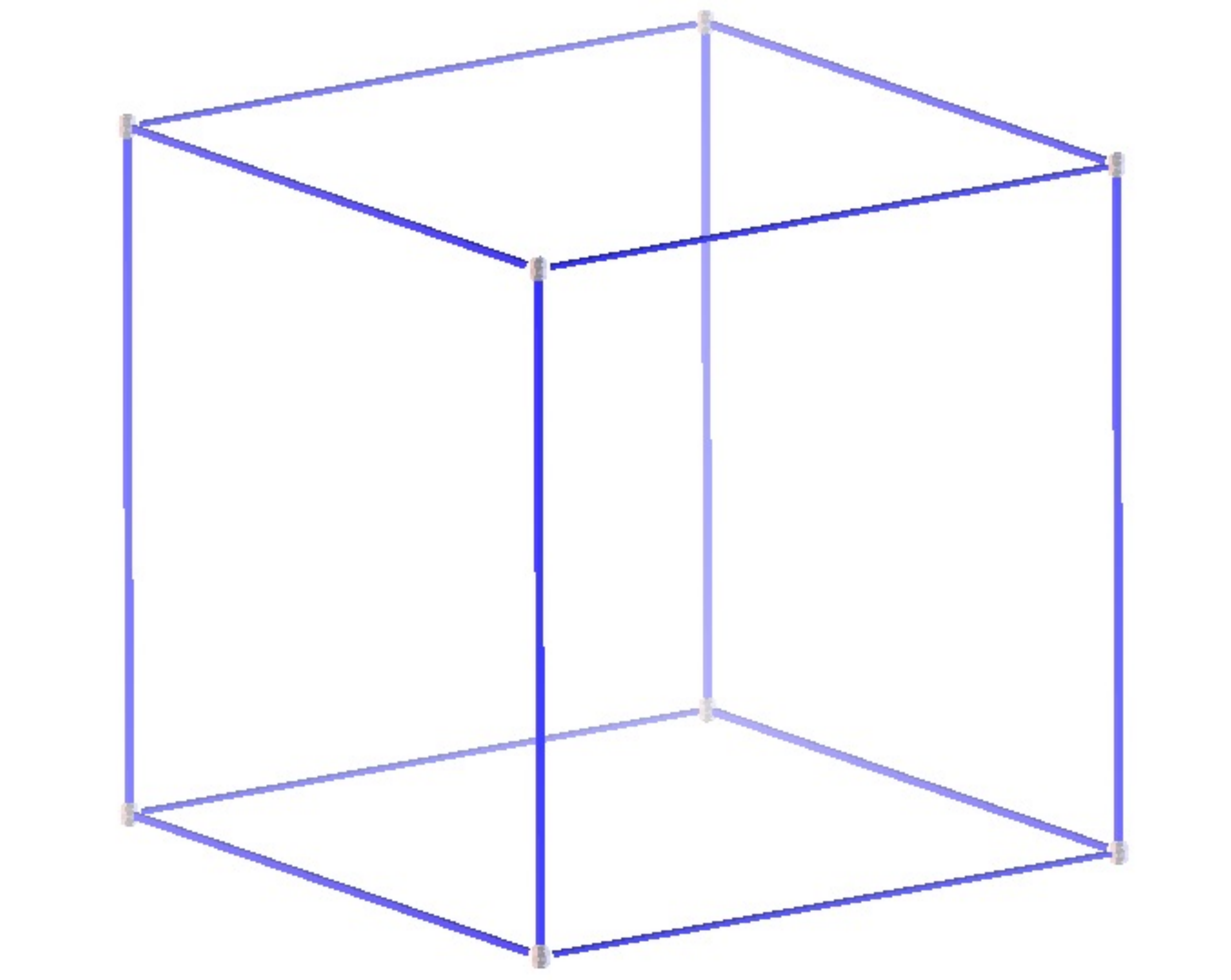}
\end{center}
\caption{The cube}\label{cubo}
\end{minipage}
\end{figure}

Let us begin with the regular tetrahedron $\Delta$ (see Figure~\ref{tetraedro}) having vertices
$$
\begin{array}{l}
\nu_1=(1,1,1)\\
\nu_2=(1,-1,-1)\\
\nu_3=(-1 ,1,-1)\\
\nu_4=(-1,-1,1).
\end{array}
$$
Consider the sublattice $L$ of $\Z^3$ that is generated by the corresponding four vectors
$$
\begin{array}{l}
Y_1=(1,1,1)\\
Y_2=(1,-1,-1)\\
Y_3=(-1 ,1,-1)\\
Y_4=(-1,-1,1).
\end{array}
$$
Notice that $Y_1+Y_2+Y_3+Y_4=0$, therefore any three of these four vectors form a basis of $L$. Moreover,
$$
\D=\bigcap_{j=1}^{4}\{\;\mu\in(\R^3)^*\;|\;\langle\mu,X_j\rangle\geq-1\;\},
$$
where $X_j=Y_j$, $j=1,\ldots,4$. Thus the regular tetrahedron $\Delta$ satisfies the hypotheses of Delzant's Theorem \cite{d}
with respect to $L$.
From the symplectic viewpoint, it is readily verified here that 
$N=\{\,(e^{2\pi i \theta},e^{2\pi i \theta},e^{2\pi i \theta},e^{2\pi i \theta})\,|\, \theta\in \R\,\}=S^1$ and therefore that
$$M_{\Delta}=\frac{\{\vz\in\C^4\,|\,|z_1|^2+|z_2|^2+|z_3|^2+|z_4|^2=4\}}{S^1}=\frac{S^7}{S^1},$$
where $S^7$ denotes the $7$--sphere of radius $2$. From the complex viewpoint, it is easy to see that
$\C^4_{\Delta}=\C^4\smallsetminus\{0\}$, that $N_{\SC}=\C^*$ and therefore that 
$$X_{\Delta}=\frac{\C^4\smallsetminus\{0\}}{\C^*}=\C\P^3.$$
Consider now the cube $\Delta$ having vertices $(\pm 1,\pm 1,\pm 1)$  (see Figure~\ref{cubo}). 
Notice that
$$
\D=\bigcap_{j=1}^{6}\{\;\mu\in(\R^3)^*\;|\;\langle\mu,X_j\rangle\geq-1\;\},
$$
where
$X_1=e_1$, $X_2=-e_1$, $X_3=e_2$, $X_4=-e_2$, $X_5=e_3$, $X_6=-e_3$. 
We can again apply the Delzant procedure, this time relatively to the lattice $\Z^3$, and we get
the $3$--dimensional group $$N=\{\,(e^{2\pi i \theta_1},e^{2\pi i \theta_1},e^{2\pi i \theta_2},e^{2\pi i \theta_2},e^{2\pi i \theta_3},e^{2\pi i \theta_3})\,|\, \theta_1,\theta_2,\theta_3\in \R\,\}=S^1\times S^1\times S^1\subset T^6$$
and, therefore, the symplectic toric manifold
$$M_{\Delta}=\frac{\{\vz\in\C^6\,|\,|z_1|^2+|z_2|^2=2,|z_3|^2+|z_4|^2=2,|z_5|^2+|z_6|^2=2\}}{N}=S^2\times S^2 \times S^2,$$
where the $S^2$'s have all radius $\sqrt{2}$. 
The corresponding complex toric manifold, on the other hand, is given by
$$X_{\Delta}=\left((\C^2\setminus\{0\})\times (\C^2\setminus\{0\})\times (\C^2\setminus\{0\})\right)/(\C^*\times\C^*\times\C^*)=
\C\P^1\times\C\P^1\times\C\P^1.$$ 
This provides an elementary example of a general fact: the symplectic toric manifold depends on the polytope, 
while the complex toric manifold only depends on the fan that is normal to the polytope. For instance, let us consider the cube
having vertices $(\pm a,\pm a,\pm a)$, with $a$ a positive real number, with the same vectors $X_1,\ldots,X_6$ as above. Then
the corresponding symplectic toric manifold is the product of three spheres of radius $\sqrt{2a}$. 
In conclusion, the symplectic structure varies, while the complex toric manifold remains the same.

\section{Nonsimple and rational: the regular octahedron}\label{octahedron}
\begin{figure}[h]
\begin{center}
\includegraphics[width=60mm]{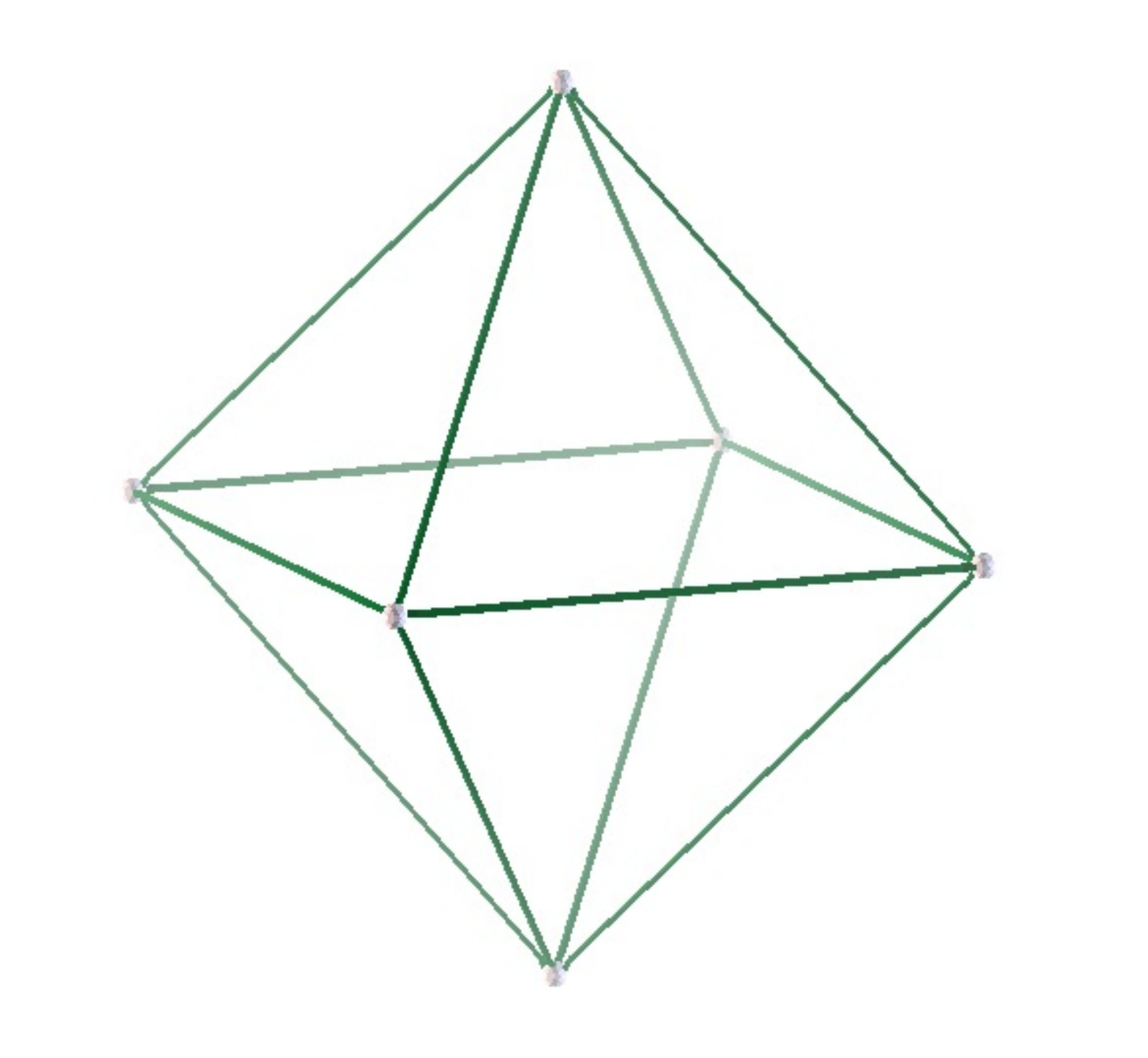}
\end{center}
\caption{The regular octahedron}\label{ottaedro}
\end{figure}

Let us consider the regular octahedron $\Delta$ that is dual to the cube of Section~\ref{smooth} (see Figure~\ref{ottaedro}). 
Its vertices are given by 
$\nu_1=e_1,\nu_2=-e_1,\nu_3=e_2,\nu_4=-e_2,\nu_5=e_3,\nu_6=-e_3$.
This polytope is rational with respect to the the lattice $L$ defined in section \ref{smooth}. 
In fact
$$
\D=\bigcap_{j=1}^{8}\{\;\mu\in(\R^3)^*\;|\;\langle\mu,X_j\rangle\geq-1\;\},
$$
where $X_i=Y_i$, $X_{4+i}=-Y_i$,  $i=1,\ldots,4$. 
However $\Delta$ is not simple.
More precisely, consider the planes 
$$H_i=\left\{\;\mu\in(\R^3)^*\;|\;\langle\mu,X_j\rangle=-1\;\right\},$$ $i=1,\ldots,8$, containing the eight facets of the octahedron.
Notice that each such plane contains exactly three vertices:
\begin{center}
\begin{tabular}{|c|c||c|c|}
\hline
plane&vertices&plane&vertices\\
\hline
$H_1$&$(\nu_2,\nu_4,\nu_6)$& $H_5$&$(\nu_1,\nu_4,\nu_5)$\\
\hline
$H_2$&$(\nu_2,\nu_3,\nu_5)$&$H_6$&$(\nu_1,\nu_4,\nu_6)$\\
\hline
$H_3$&$(\nu_1,\nu_4,\nu_5)$&$H_7$&$(\nu_2,\nu_3,\nu_6)$\\
\hline
$H_4$&$(\nu_1,\nu_3,\nu_6)$&$H_8$&$(\nu_2,\nu_4,\nu_5)$\\
\hline
\end{tabular}
\end{center}
while each vertex is given by the intersection of four planes:
\begin{equation}
\label{carteottaedro}
\begin{array}{c}
\nu_1=H_3\cap H_4\cap H_5\cap H_6\\
\nu_2=H_1\cap H_2\cap H_7\cap H_8\\
\nu_3=H_2\cap H_4\cap H_5\cap H_7\\
\nu_4=H_1\cap H_3\cap H_6\cap H_8\\
\nu_5=H_2\cap H_3\cap H_5\cap H_8\\
\nu_6=H_1\cap H_4\cap H_6\cap H_7.
\end{array}
\end{equation}
Therefore, each vertex is a singular face of $\Delta$.
If we apply the generalized Delzant construction we find that $N=\exp(\n)$, with
$$\n=\{(\theta_1,\theta_2,\theta_3,\theta_4,\theta_1+s,\theta_2+s,\theta_3+s,\theta_4+s)\;|\;\theta_1,\theta_2,\theta_3,\theta_4,s\in \R\}.$$
The level set $\Psi^{-1}(0)$ for the moment mapping  $\Psi$ with respect to the induced $N$--action 
on $\C^8$ is therefore given by the compact subset of points $\underline{z}$ in $\C^8$ such that
$$
\left\{
\begin{array}{l}
|z_1|^2+|z_5|^2=2\\
|z_2|^2+|z_6|^2=2\\
|z_3|^2+|z_7|^2=2\\
|z_4|^2+|z_8|^2=2\\
|z_3|^2+|z_4|^2=|z_5|^2+|z_6|^2.
\end{array}
\right.
$$
Notice that the $N$--action has isotropies of dimension $1$ at the
$T^8$--orbits corresponding to the $6$ (singular) vertices. For example, consider the first vertex: $\nu_1=H_3\cap H_4\cap H_5\cap H_6$.
The corresponding orbit in $\Psi^{-1}(0)$ is
$$\{\underline{z}\in\C^8\;|\;|z_1|^2=2,|z_2|^2=2,|z_7|^2=2,|z_8|^2=2\},$$
with isotropy
$\{(1,1,e^{2\pi i t},e^{2\pi i t},1,1,e^{2\pi i t},e^{2\pi i t})\;|\; t\in \R\}$.
Away from the orbits corresponding to the vertices, the level set is smooth.
By Theorem~\ref{nonsemplice-re}, the orbit space $$M_{\Delta}=\Psi^{-1}(0)/N$$ is an elementary example of symplectic toric space stratified by manifolds. The six orbits in $\Psi^{-1}(0)$ corresponding to the vertices yield six singular points in the quotient $M$.
Their complement is the regular stratum of $M$, which is a smooth symplectic manifold.

From the complex viewpoint, by (\ref{aperto}), the set
$\C^8_{\Delta}$ is the union of $6$ open sets, each of which corresponds to a vertex. The Cox
quotient $X_{\Delta}=\C^8_{\Delta}//N_{\SC}$
is an elementary example of complex toric space stratified by manifolds. 
The quotients $M_{\Delta}$ and $X_{\Delta}$ can be identified according to Theorem~\ref{non semplice iso}.

Let us briefly describe $M_{\Delta}$ as a symplectic toric stratified space. By Theorem~\ref{nonsemplice-re}, the action of the torus $D^3=\R^3/L$ on $M_{\Delta}$ is effective and continuous. Moreover,
the image of the continuous mapping $\Phi=(\pi^*)^{-1}\circ J$ is exactly $\Delta$ and 
each singular point is sent by $\Phi$ to the corresponding vertex. The action of the torus $D^3=\R^3/L$
on the regular stratum is smooth and Hamiltonian, and the mapping $\Phi$ is the moment mapping with respect to the $D^3$--action. By Theorem~\ref{non semplice iso}, the regular stratum is a $6$--dimensional K\"ahler manifold \cite[Remark~6.6]{stratificato-re}.

We now describe the local structure of the regular stratum and of the stratification around the singular points.
Consider, for example, the vertex $\nu_1$. We recall that the singular point $m_1\in M_{\Delta}$ corresponding to $\nu_1$ is given by
$$\{\underline{z}\in\C^8\;|\;|z_1|^2=2,|z_2|^2=2,|z_7|^2=2,|z_8|^2=2\}/N.$$
We compute a symplectic chart for the regular stratum in a neighborhood of $m_1$.
Consider the open subset
$$U=\{(u_3,u_4,u_5)\in\C^3\setminus\{0\}\;|\;|u_5|^2<|u_3|^2+|u_4|^2<2+|u_5|^2,|u_i|^2<2, i=3,4,5\}$$
and the mapping
$$\begin{array}{ccccc}
\tau&\colon&U&\rightarrow&
\{\underline{z}\in\Psi^{-1}(0)\;|\;z_i\neq0, i=1,2,6,7,8\;\hbox{and}\;(z_3,z_4,z_5)\neq(0,0,0)\}\\
&&(u_3,u_4,u_5)&\mapsto&\left(\tau_1(\underline{u}),
\tau_2(\underline{u}),u_3,u_4,u_5,\tau_6(\underline{u}),\tau_7(\underline{u}),\tau_8(\underline{u})\right)
\end{array},
$$
where 
$$
\begin{array}{l}
\tau_{1}(\underline{u})=\sqrt{2-|u_5|^2}\\
\tau_{2}(\underline{u})=\sqrt{2-|u_3|^2-|u_4|^2+|u_5|^2}\\
\tau_{6}(\underline{u})=\sqrt{|u_3|^2+|u_4|^2-|u_5|^2}\\
\tau_{7}(\underline{u})=\sqrt{2-|u_3|^2}\\
\tau_{8}(\underline{u})=\sqrt{2-|u_4|^2}.\\
\end{array}
$$
The mapping $\tau$ is a homeomorphism (see Proof of \cite[Theorem~5.3]{stratificato-re}).
The pair $((U,\omega_{0}),\tau)$, where $\omega_{0}$ is the standard symplectic structure of $\C^3$, defines a symplectic chart for the regular part in a neighborhood of $m_1$. To find other charts, repeat the procedure by replacing the indices $\{3,4,5\}$ with any
triple of indices $\{i,j,k\}$ contained in one of the six sets 
\begin{equation}\label{sets}
\{3,4,5,6\},\{1,2,7,8\},\{2,4,5,7\},\{1,3,6,8\},\{2,3,5,8\},\{1,4,6,7\}.
\end{equation}
We recall that, according to (\ref{carteottaedro}), each of these sets correspond  to a different vertex. 
We need at least $12$ charts to obtain a symplectic atlas for the regular stratum.

The regular part can also be seen as a complex manifold. Consider the open subset
$$V=\C^3\smallsetminus\{0\}$$
and the mapping
$$\begin{array}{ccccc}
\tau_{\SC}&\colon&V&\rightarrow&
\{\underline{z}\in\C^8_{\Delta}\;|\;z_i\neq0, i=1,2,6,7,8\;\hbox{and}\;(z_3,z_4,z_5)\neq(0,0,0)\}\\
&&(u_3,u_4,u_5)&\mapsto&\left(1,
1,u_3,u_4,u_5,1,1,1\right)
\end{array}.
$$
The mapping $\tau_{\SC}$ is a homeomorphism.
The pair $(V,\tau_{\SC})$ defines a complex chart for the regular part in a neighborhood of  the singular point $x_1=\chi(m_1)$.
Again, a complex atlas can be obtained  
by considering all the triples of indices $\{i,j,k\}$ contained in one
of the six sets (\ref{sets}).
Analogously to \cite[Example~3.8]{cx}, one can compute the local expression of the mapping $\chi$ as a diffeomorphism from $U$ to $V$ and find the local expression for the K\"ahler form on $V$.

Let us now give an example of the existence of non--closed $A$--orbits for nonsimple polytopes.
Consider a point in $\C^8_{\Delta}$ that does not lie in any $\C^8_{F}$, 
for instance $$\vz=(z_1:z_2:z_3:0:0:0:z_7:z_8),\quad z_1,z_2,z_3,z_7,z_8\neq0.$$
Take the element $a(t)=(1,1,e^{-2\pi t},e^{-2\pi t},e^{2\pi t},e^{2\pi t},1,1)\in A$. When $t$ tends to $+\infty$,
$a(t)\cdot\vz$ tends to $(z_1:z_2:0:0:0:0:z_7:z_8)$, which therefore lies in the closure of the orbit $A\cdot\vz$. Notice
that, by Theorem~\ref{closedorbits}, the orbit through $(z_1:z_2:0:0:0:0:z_7:z_8)$ is the only closed orbit contained in
$\overline{A\cdot\vz}$.

Let us now show that the quotient $M_{\Delta}\simeq X_{\Delta}$ is a stratified space in the sense of \cite{gmp} with isolated singularities. More precisely we show that, in this case, a neighborhood of each singular point $m_j$, with $j=1,\ldots,6$, can be identified with a cone over a manifold, called the link of $m$.
We consider the first vertex; the same argument applies to the others. Consider the cone stemming from $\nu_1$
$$
{\cal C}=\bigcap_{j=3}^{6}\{\;\mu\in(\R^3)^*\;|\;\langle\mu,X_j\rangle\geq-1\;\}.
$$
The symplectic and complex toric spaces corresponding to this cone are
$X_{\cal C}=\C^4//N_{\SC}({\cal C})$ 
and
$M_{\cal C}=\{\vz\in\C^4\;|\;|z_3|^2+|z_4|^2=|z_5|^2+|z_6|^2\}/N({\cal C}),$
where $$N({\cal C})=\exp\{(s,s,-s,-s)\;|\;s\in\R\}.$$
We now see explicitly that $M_{\cal C}$ and $X_{\cal C}$ are cones (diffeomorphic by the proof of \cite[Theorem~3.3]{stratificato-cx}) which are local models for our toric space
$M_{\Delta}\simeq X_{\Delta}$ near the singular point $m_1=\chi^{-1}(x_1)$. This will prove that $M_{\Delta}\simeq X_{\Delta}$
is a stratified space. Following the recipes given in \cite[Sections~3.2,3.4]{stratificato-cx},
the link for $m_1=\chi^{-1}(x_1)$ can be found by cutting the cone $\cal C$ with 
the $yz$--plane; this gives the square $\Delta_L$ of vertices $(0,\pm 1,\pm 1)$.
The $5$--dimensional {\em real link} is then given by
$$L=\begin{array}{l}\{\vz\in\C^4\;|\;|z_3|^2+|z_4|^2=|z_5|^2+|z_6|^2=2\}/N({\cal C})
\simeq (S^3\times S^3)/S^1,\end{array}$$
where the spheres $S^3$ have both radius $\sqrt{2}$ and
$S^1$ acts with weights $+1$ on the first two coordinates and with weights $-1$ on the last two.
From the complex point of view, we have
$$
L_{\SC}=
\left((\C^2\smallsetminus\{0\})\times (\C^2\smallsetminus\{0\})\right)/(N({\cal C})\exp\{i(s,s,s,s)\;|\;s\in\R\}).
$$
Notice that the inclusion $\{\vz\in\C^4\;|\;|z_3|^2+|z_4|^2=|z_5|^2+|z_6|^2=2\}\hookrightarrow\C^4$ induces a diffeomorphism
$L\simeq L_{\SC}$.
By applying the Delzant procedure to the square $\Delta_L$, with respect to the the normal vectors $\pm e_2,\pm e_3$,
one finds $N(\Delta_L)=\exp\{(s,s,t,t)\;|\;s,t\in\R\}$. Thus the symplectic toric manifold $M_{\Delta_L}$ is given by
$$\{\vz\in\C^4\;|\;|z_3|^2+|z_4|^2=|z_5|^2+|z_6|^2=2\}/N({\Delta_L})\simeq S^2\times S^2,$$
where the spheres $S^2$ also have radius $\sqrt{2}$.
We follow the terminology of \cite{sl} and refer to this manifold as the {\em symplectic link}.
Correspondingly, {\em the complex link} in the sense of \cite[p.15]{gmp}, is given by
$$X_{\Delta_L}=\left((\C^2\smallsetminus\{0\})\times (\C^2\smallsetminus\{0\})\right)/N_{\SC}(\Delta_L)\simeq\C\P ^1\times\C\P^1.$$
Notice that 
$$\C^4//N_{\SC}({\cal C})\smallsetminus[0]=\left((\C^2\smallsetminus\{0\})\times (\C^2\smallsetminus\{0\})\right)/N_{\SC}({\cal C});$$
this shows that the cone $X_{\cal C}$ is a complex cone over the complex link.
Finally, remark that our space $M_{\cal C}\simeq X_{\cal C}$, without the singular point that corresponds to the cone apex,
projects first onto the real link $L\simeq L_{\SC}$ and then onto the symplectic and complex links as follows:
\begin{equation}\label{diagramma reale}
\xymatrix{
\{\vz\in\C^4\smallsetminus\{0\}\;|\;|z_3|^2+|z_4|^2=|z_5|^2+|z_6|^2\}/N({\cal C})
\ar@{<->}[r]
\ar[d]^{p^1}
\ar@/_3pc/[dd]_p
&
\C^4//N_{\SC}({\cal C})\smallsetminus[0]
\ar[d]^{p^1_{\SC}}
\ar@/^3pc/[dd]^{p_{\SC}}\\
L
\ar@{<->}[r]
\ar[d]^{p^2}&
L_{\SC}
\ar[d]^{p^2_{\SC}}\\
S^2\times S^2\ar@{<->}[r]&\C\P ^1\times\C\P^1
}
\end{equation}
On the left hand side,
$p^1([\vz])=[rz_1:rz_2:rz_3:rz_4]$, 
with $r=\left(\frac{1}{2}(|z_3|^2+|z_4|^2)\right)^{-1/2}$, $p^2$ is the natural projection and $p$ is the composite. 
On the right hand side the mappings are the natural projections.

The fiber of the projections $p^1$ and $p^1_{\SC}$ is $\R_{>0}$,
the fiber of the projections $p^2$ and $p^2_{\SC}$ is $$N(\Delta_L)/N({\cal C})\simeq \exp\{(t,t,t,t)\;|\;t\in\R\}\simeq S^1,$$
and the fiber of the projections $p$ and $p_{\SC}$ is
$$N_{\SC}({\Delta_L})/N_{\SC}({\cal C})\simeq\exp\{(u,u,u,u)\;|\;u\in\C\}\simeq\R_{>0}\times S^1\simeq\C^*.$$
Therefore, $X_{\cal C}\simeq M_{\cal C}$ is a real cone over $L_{\SC}\simeq L$ and a
complex cone over the symplectic/complex link. Finally, notice that $L/S^1$ gives the symplectic link.

We now show that  $M_{\cal C}$ and $X_{\cal C}$ are local models for $M_{\Delta}$ and $X_{\Delta}$ respectively.
Consider, in the cone $\{\vz\in\C^4\smallsetminus\{0\}\;|\;|z_3|^2+|z_4|^2=|z_5|^2+|z_6|^2\}/N({\cal C})$, the neighborhood ${\cal C}_1$ of
the apex given by the set of points such that $|z_i|^2<2, i=3,4,5,6$.
The continuos mapping
$$
\begin{array}{ccc}
{\cal C}_1&\longrightarrow&
\{\vz\in\Psi^{-1}(0)\;|\;z_i\neq0, i=1,2,7,8\}\\
{[}z_3:z_4:z_5:z_6]&\longmapsto&[\sqrt{2-|z_3|^2}:\sqrt{2-|z_4|^2}:z_3:z_4:z_5:z_6:\sqrt{2-|z_7|^2}:\sqrt{2-|z_8|^2}]
\end{array}$$
is a symplectomorphism on the regular part and sends singular point to singular point. 
On the other hand, the continuos mapping
$$
\begin{array}{ccc}
\C^4//N^{\cal C}_{\SC}&\longrightarrow&\{\vz\in\C^8\;|\;z_i\neq0, i=1,2,7,8\}//N_{\SC}\subset \C^8_{\Delta}//N_{\SC}\\
{[}z_3:z_4:z_5:z_6]&\longmapsto&[1:1:z_3:z_4:z_5:z_6:1:1]
\end{array}
$$
is a biholomorphism on the regular part and also sends singular point to singular point. 

\section{Simple and nonrational: the regular dodecahedron}\label{dodecahedron}
\begin{figure}[h]
\begin{center}
\includegraphics[width=70mm]{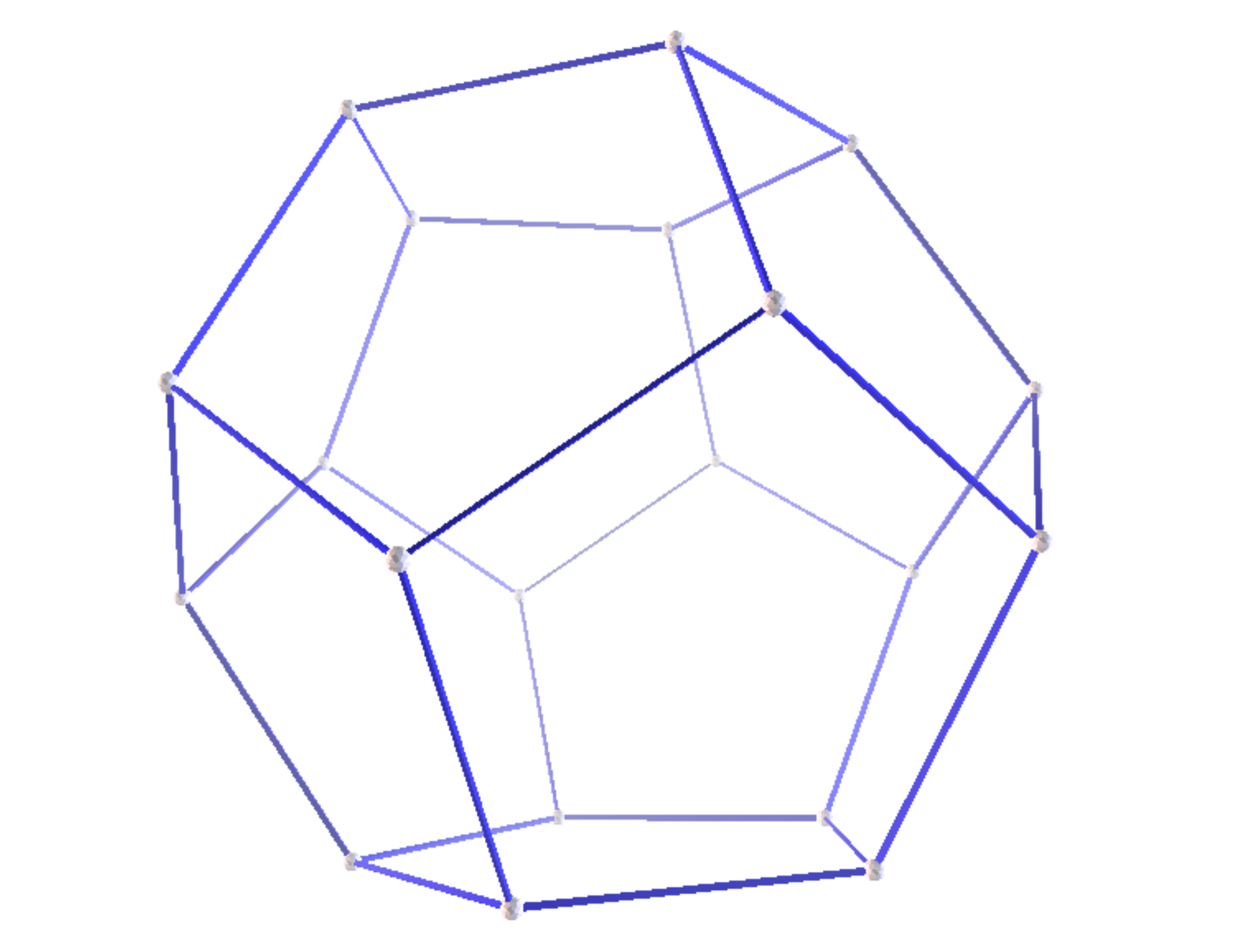}
\end{center}
\caption{The regular dodecahedron}\label{dodecaedro}
\end{figure}

\begin{figure}[h]
\begin{minipage}[b]{0.45\linewidth}
\begin{center}
\includegraphics[width=60mm]{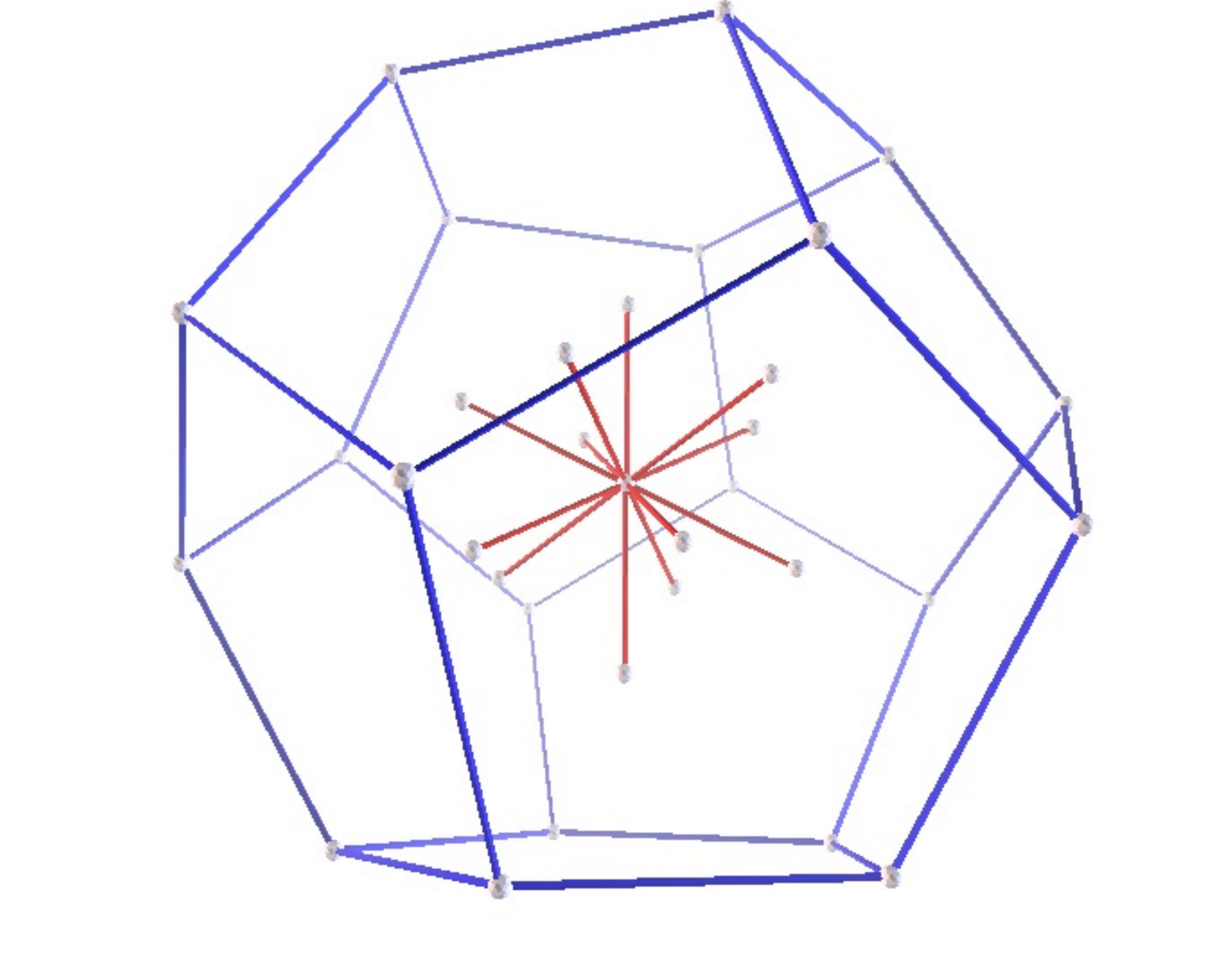}
\end{center}
\caption{The vectors $X_1,\ldots,X_{12}$}\label{stelladirossi}
\end{minipage}
\qquad
\begin{minipage}[b]{0.45\linewidth}
\begin{center}
\includegraphics[width=60mm]{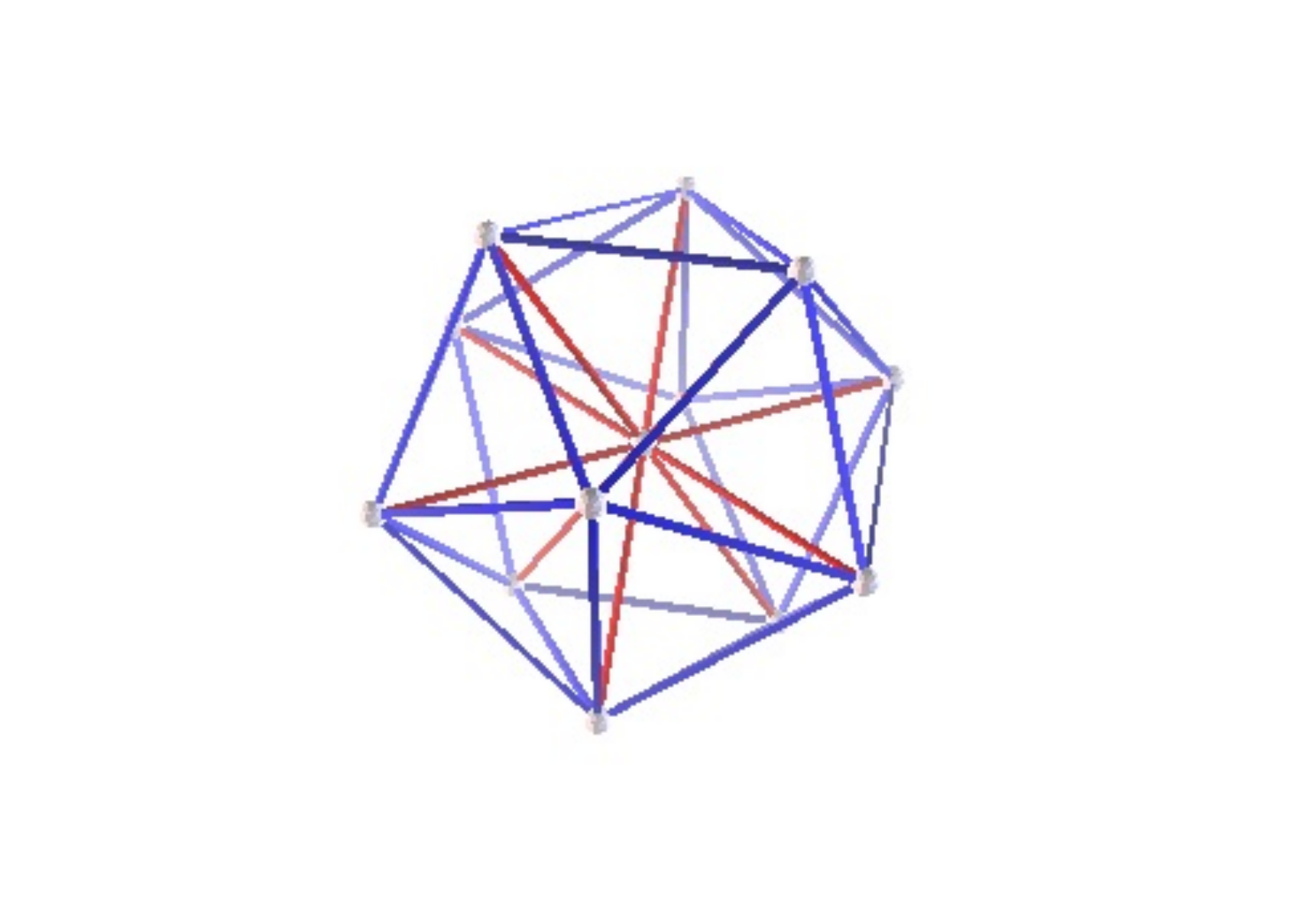}
\end{center}
\caption{The dual icosahedron}\label{dualsenzagialli}\label{stargrid}
\end{minipage}
\end{figure}

Let $\phi=\frac{1+\sqrt{5}}{2}$ be the {\em golden ratio} and remark that it satisfies the equation $\phi=1+\frac{1}{\phi}$.
Let $\Delta$ be the regular dodecahedron  having vertices
$$
\begin{array}{l}
(\pm 1,\pm 1,\pm 1)\\
(0,\pm \phi,\pm \frac{1}{\phi})\\
(\pm \frac{1}{\phi} ,0,\pm \phi)\\
(\pm \phi,\pm \frac{1}{\phi},0)
\end{array}
$$
(see Figure~\ref{dodecaedro}).
The polytope $\Delta$ is simple but not rational. 
However, it is quasirational with respect to the quasilattice $P$,
known in physics as the {\em simple icosahedral lattice} \cite{rmw},
that is generated by the following vectors in $\R^3$:
$$
\begin{array}{l}
Y_1=(\frac{1}{\phi^2},\frac{1}{\phi},0)\\
Y_2=(0,\frac{1}{\phi^2},\frac{1}{\phi})\\
Y_3=(\frac{1}{\phi},0,\frac{1}{\phi^2})\\
Y_4=(-\frac{1}{\phi^2},\frac{1}{\phi},0)\\
Y_5=(0,-\frac{1}{\phi^2},\frac{1}{\phi})\\
Y_6=(\frac{1}{\phi},0,-\frac{1}{\phi^2}).
\end{array}
$$
In fact, a straightforward computation shows that:
$$
\D=\bigcap_{j=1}^{12}\{\;\mu\in(\R^3)^*\;|\;\langle\mu,X_j\rangle\geq-1\;\},
$$
where $X_i=Y_i$, $X_{6+i}=-Y_i$, $i=1,\ldots,6$. 
Notice that the vectors $X_1,\ldots,X_{12}$ point to the twelve vertices of the regular icosahedron that is dual to $\D$ 
(see Figures~\ref{stelladirossi} and \ref{dualsenzagialli}).
We recall from \cite{pdoc} that, if we apply the generalized Delzant construction to $\Delta$ with respect to the quasilattice $P$ and the vectors $X_1,\ldots,X_{12}$,
we get $N=\exp(\n)$, $\n$ being the $9$--dimensional subspace of $\R^{12}$ that is spanned by the vectors 
$$
\begin{array}{l}
e_1+e_7\\
e_2+e_8\\
e_3+e_9\\
e_4+e_{10}\\
e_5+e_{11}\\
e_6+e_{12}\\e_1+e_2-\phi(e_3+e_4)\\
e_2+e_3-\phi(e_1+e_5)\\
e_1+e_3-\phi(e_2+e_6).
\end{array}
$$
The level set $\Psi^{-1}(0)$ for the moment mapping  $\Psi$ with respect to the induced $N$--action on $\C^{12}$ is therefore given by the compact
subset of points $\underline{z}$ in $\C^{12}$ such that
$$
\left\{
\begin{array}{lll}
|z_1|^2+|z_7|^2&=&2\\
|z_2|^2+|z_8|^2&=&2\\
|z_3|^2+|z_9|^2&=&2\\
|z_4|^2+|z_{10}|^2&=&2\\
|z_5|^2+|z_{11}|^2&=&2\\
|z_6|^2+|z_{12}|^2&=&2\\
|z_1|^2+|z_2|^2&-&\phi(|z_3|^2+|z_4|^2)=-\frac{2}{\phi}\\
|z_2|^2+|z_3|^2&-&\phi(|z_1|^2+|z_5|^2)=-\frac{2}{\phi}\\
|z_1|^2+|z_3|^2&-&\phi(|z_2|^2+|z_6|^2)=-\frac{2}{\phi}.
\end{array}
\right.
$$  
The quotient $M_{\Delta}=\Psi^{-1}(0)/N$ is a symplectic quasifold; it has an atlas made of $20$ charts, each corresponding to a different fixed point of the $D^3$--action; we refer the reader to \cite{pdoc} for a description of one of them.

From the complex viewpoint, the complex toric quasifold corresponding to the dodecahedron, with the choice of
the same vectors $X_j$ ($j=1,\ldots,12$) above is given by the quotient
$$X_{\Delta}=\C^{12}_{\Delta}/N_{\SC},$$ where $\C^{12}_{\Delta}$ is the open subset of $\C^{12}$ given by the union of the $20$ 
open subsets defined in (\ref{aperto}). The symplectic and complex quotients can be identified by Theorem~\ref{teoremadellachi}.

Let us describe a chart for the complex toric quasifold $X_{\Delta}$ around the fixed point corresponding to the vertex $(-1,-1,-1)$.
Let ${\tilde V}=\C^3$ and consider the following slice of $\C^{12}_{\Delta}$ that is transversal to the $N_{\SC}$--orbits
$$\begin{array}{ccccc}
\tilde{\tau}_{\SC}&\colon&\tilde{V}&\rightarrow&
\{\vw\in\C^{12}_{\Delta}\;|\;w_i\neq0, i=4,\ldots ,12\}\\
&&(z_1,z_2,z_3)&\mapsto&\left(z_1,z_2,z_3,1,\ldots,1\right)
\end{array},
$$
The mapping $\tilde{\tau}_{\SC}$ induces a homeomorphism 
$$
\begin{array}{ccccc}
\tau_{\SC}&\colon&\tilde{V}/\Gamma&\longrightarrow& U\\
&&\,[(z_1,z_2,z_3)]&\longmapsto&[{\tilde{\tau}_{\SC}}(z_1,z_2,z_3)]
\end{array},
$$
where the open subset $V$ of $X_{\Delta}$ is the quotient
$$\{\vw\in\C^{12}_{\Delta}\;|\;w_i\neq0, i=4,\ldots ,12\}/N_{\SC}$$
and the discrete group $\Gamma$ is given by
$$
\Gamma=\left\{\,(e^{2\pi i\phi(h+l)},e^{2\pi i\phi(h+k)},e^{2\pi i\phi(k+l)})\in T^3\;|\; h,k,l\in\Z\right\}.
$$
The triple $(V,\tau,\tilde{V}/\Gamma)$ defines a complex chart for $X_{\Delta}$. The others can be described similarly.

\section{Nonsimple and nonrational: the regular icosahedron}\label{icosahedron}
\begin{figure}[h]
\begin{center}
\includegraphics[height=45mm]{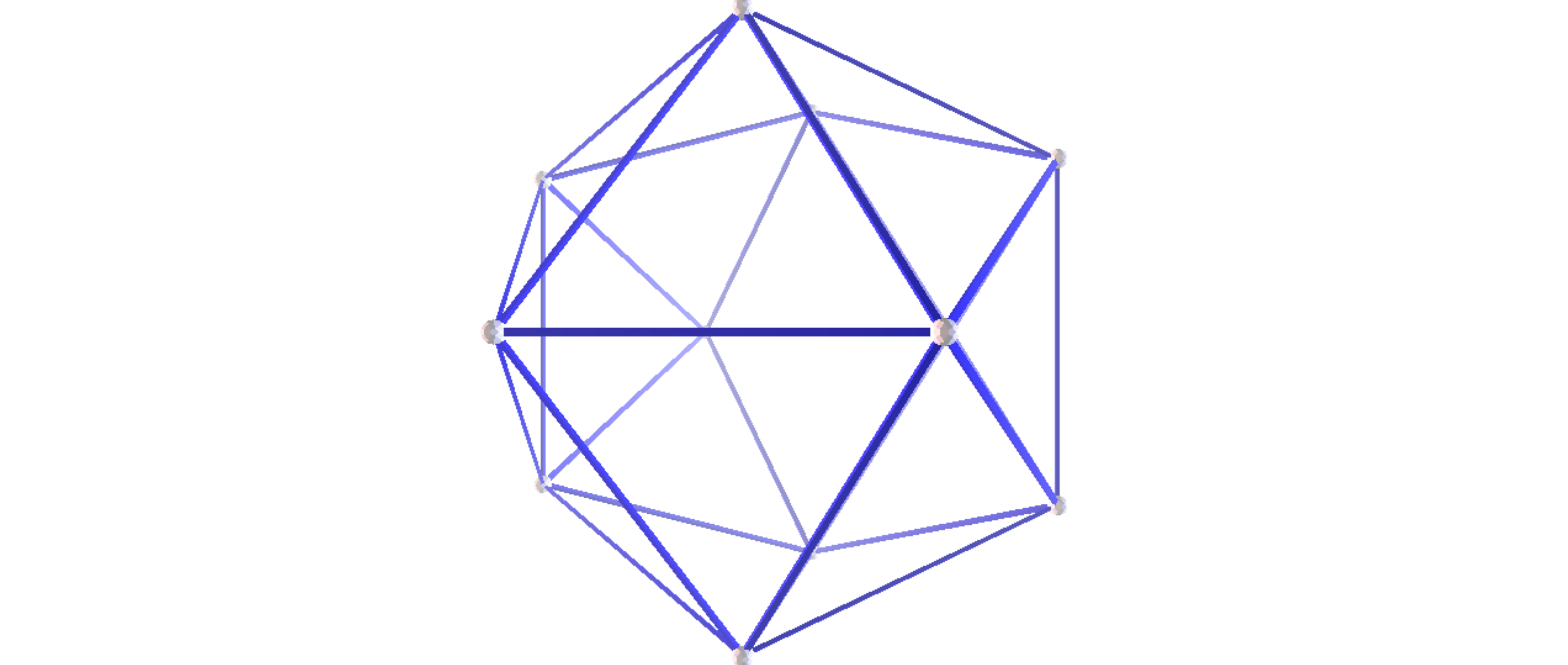}
\end{center}
\caption{The regular icosahedron}\label{icosaedro}
\end{figure}

\begin{figure}[h]
\begin{minipage}[b]{0.45\linewidth}
\begin{center}
\includegraphics[width=60mm]{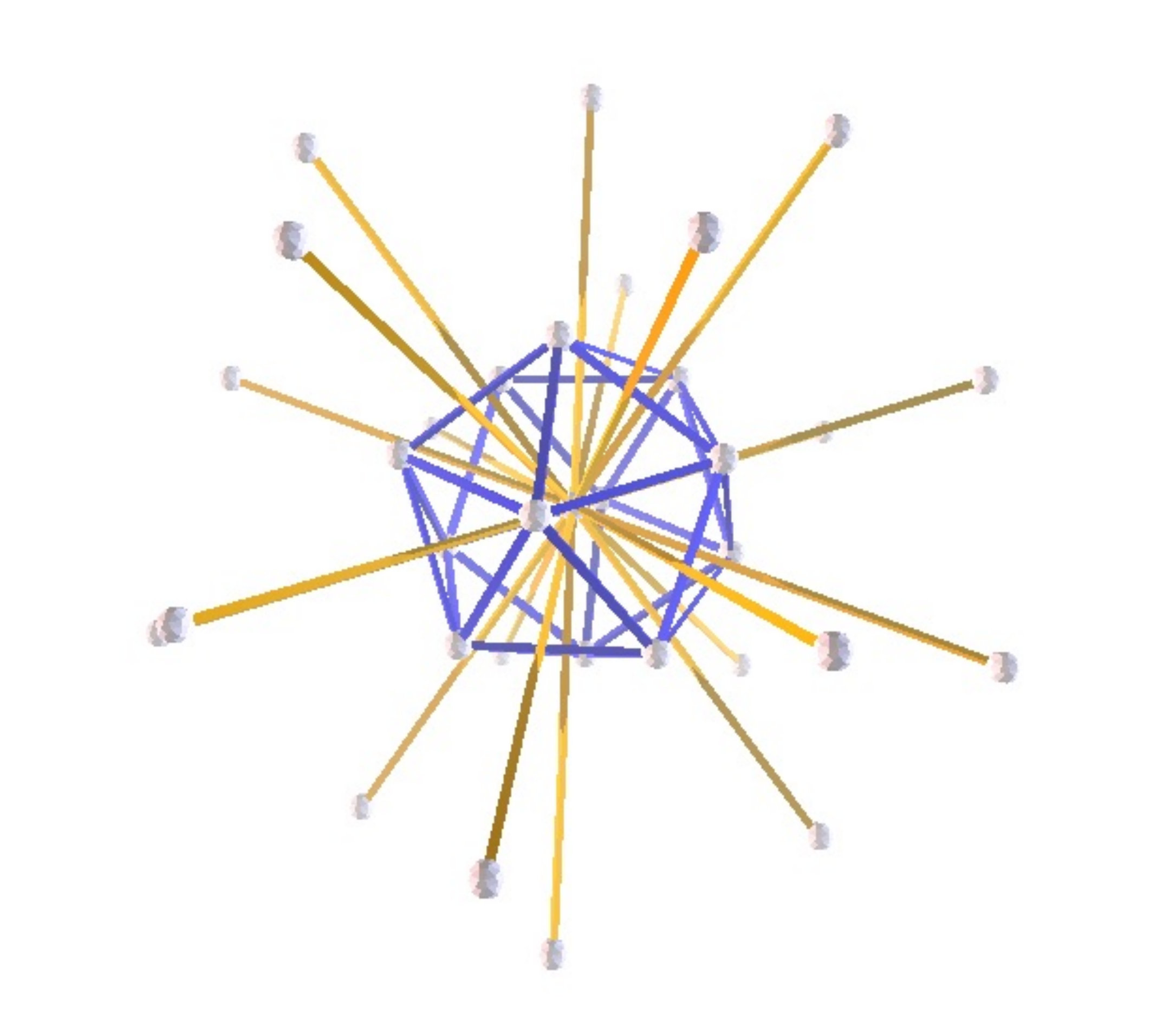}
\end{center}
\caption{The normal vectors $X_1,\ldots,X_{20}$}\label{stelladigialli}
\end{minipage}
\qquad
\begin{minipage}[b]{0.45\linewidth}
\begin{center}
\includegraphics[width=60mm]{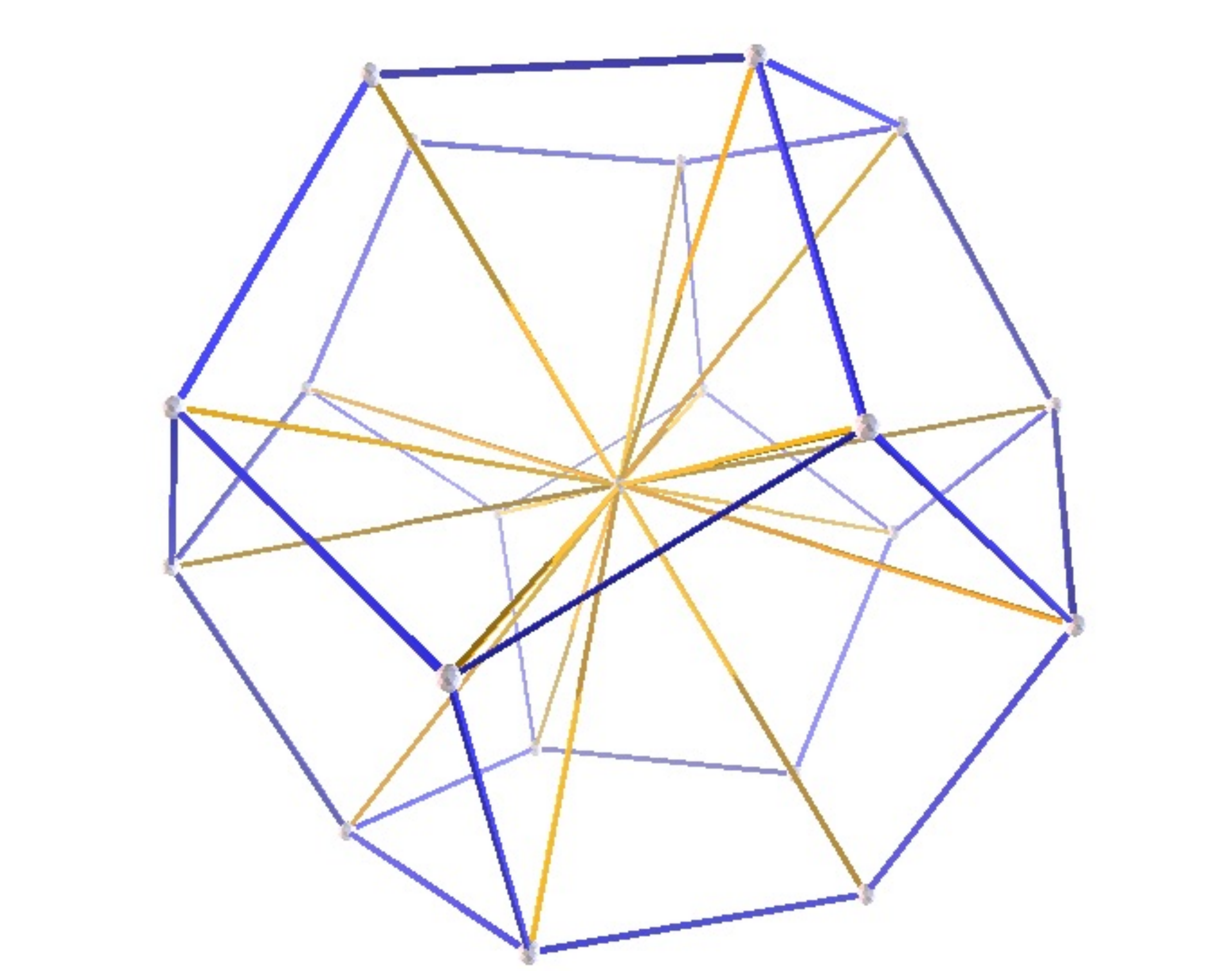}
\end{center}
\caption{The dual dodecahedron}\label{dualsenzarossi}
\end{minipage}
\end{figure}

Let $\Delta$ be the regular icosahedron that is dual to the dodecahedron in Section \ref{dodecahedron} (see Figure~\ref{icosaedro}). 
Its vertices are given by
$$
\begin{array}{ll}
\nu_1=(\frac{1}{\phi^2},\frac{1}{\phi},0)&\quad\quad\nu_7=-\nu_1\\
\nu_2=(0,\frac{1}{\phi^2},\frac{1}{\phi})&\quad\quad\nu_8=-\nu_2\\
\nu_3=(\frac{1}{\phi},0,\frac{1}{\phi^2})&\quad\quad\nu_9=-\nu_3\\
\nu_4=(-\frac{1}{\phi^2},\frac{1}{\phi},0)&\quad\quad\nu_{10}=-\nu_4\\
\nu_5=(0,-\frac{1}{\phi^2},\frac{1}{\phi})&\quad\quad\nu_{11}=-\nu_5\\
\nu_6=(\frac{1}{\phi},0,-\frac{1}{\phi^2})&\quad\quad\nu_{12}=-\nu_6.
\end{array}
$$
It is not a rational polytope. However it is quasirational with respect to the quasilattice $B$,
known in physics as the {\em body--centered lattice} \cite{rmw}, that is generated by the six vectors
$$
\begin{array}{l}
V_1=(\phi,\frac{1}{\phi},0)\\
V_2=(0,\phi,\frac{1}{\phi})\\
V_3=(\frac{1}{\phi},0,\phi)\\
V_4=(-\phi,\frac{1}{\phi},0)\\
V_5=(0,-\phi,\frac{1}{\phi})\\
V_6=(\frac{1}{\phi},0,-\phi).
\end{array}
$$
Consider in fact the four additional vectors in $B$ given by
$$
\begin{array}{lllll}
V_7&=&-V_4-V_5-V_6&=&(1,1,1)\\
V_8&=&-V_1+V_3-V_5&=&(-1,1,1)\\
V_9&=&+V_1-V_2-V_6&=&(1,-1,1)\\
V_{10}&=&+V_2-V_3-V_4&=&(1,1,-1).\\
\end{array}
$$
Then
$$
\D=\bigcap_{j=1}^{20}\left\{\;\mu\in(\R^3)^*\;|\;\langle\mu,X_j\rangle\geq-1\;\right\},
$$
where $X_i=V_i$, $X_{10+i}=-V_i$, $i=1,\ldots,10$.
Notice that the vectors $X_1,\ldots, X_{20}$ point to the vertices of the dual regular dodecahedron of Section \ref{dodecahedron} (see Figures~\ref{stelladigialli} and \ref{dualsenzarossi}).
The regular icosahedron is not a simple polytope.
More precisely, consider the planes $H_i=\left\{\;\mu\in(\R^3)^*\;|\;\langle\mu,X_j\rangle=-1\;\right\}$, $i=1,\ldots,20$, supporting the facets of $\Delta$.
Notice that each such plane contains exactly three vertices:
\begin{center}
\begin{tabular}{|c|c||c|c|}
\hline
plane&vertices&plane&vertices\\
\hline
$H_1$&$(\nu_7,\nu_9,\nu_{12})$& $H_{11}$ &$(\nu_1,\nu_3,\nu_6)$\\
\hline
$H_2$&$(\nu_7,\nu_8,\nu_{10})$& $H_{12}$ &$(\nu_1,\nu_2,\nu_4)$\\
\hline
$H_3$&$(\nu_8,\nu_9,\nu_{11})$& $H_{13}$ &$(\nu_2,\nu_3,\nu_5)$\\
\hline
$H_4$&$(\nu_3,\nu_6,\nu_{10})$& $H_{14}$ &$(\nu_4,\nu_9,\nu_{12})$\\
\hline
$H_5$&$(\nu_1,\nu_4,\nu_{11})$& $H_{15}$ &$(\nu_5,\nu_7,\nu_{10})$\\
\hline
$H_6$&$(\nu_2,\nu_5,\nu_{12})$& $H_{16}$ &$(\nu_6,\nu_8,\nu_{11})$\\
\hline
$H_7$&$(\nu_7,\nu_8,\nu_9)$& $H_{17}$ &$(\nu_1,\nu_2,\nu_3)$\\
\hline
$H_8$&$(\nu_6,\nu_8,\nu_{10})$& $H_{18}$ &$(\nu_2,\nu_4,\nu_{12})$\\
\hline
$H_9$&$(\nu_4,\nu_9,\nu_{11})$& $H_{19}$ &$(\nu_3,\nu_5,\nu_{10})$\\
\hline
$H_{10}$&$(\nu_5,\nu_7,\nu_{12})$& $H_{20}$ &$(\nu_1,\nu_6,\nu_{11})$\\
\hline
\end{tabular}
\end{center}
and that each vertex is given by the intersection of five planes:
\begin{equation}
\label{tabella vertici icosaedro}\begin{array}{l}
\nu_1=H_{5}\cap H_{11}\cap H_{12}\cap H_{17}\cap H_{20} \\
\nu_2=H_6\cap H_{12}\cap H_{13}\cap  H_{17}\cap H_{18}\\
\nu_3=H_4\cap H_{11}\cap H_{13}\cap H_{17}\cap H_{19}\\
\nu_4=H_5\cap H_{9}\cap H_{12}\cap H_{14}\cap H_{18}\\
\nu_5=H_6\cap H_{10}\cap H_{13}\cap H_{15}\cap H_{19}\\
\nu_6= H_4\cap H_8 \cap H_{11}\cap H_{16}\cap  H_{20}\\
\nu_7=H_{1}\cap H_{2}\cap H_{7}\cap H_{10}\cap  H_{15}\\
\nu_8=H_{2}\cap H_{3}\cap H_{7}\cap H_{8}\cap H_{16}\\
\nu_9=H_{1}\cap H_{3}\cap H_{7}\cap H_{9}\cap H_{14}\\
\nu_{10}=H_{2}\cap H_{4}\cap H_{8}\cap H_{15}\cap H_{19}\\
\nu_{11}=H_{3}\cap H_{5}\cap H_{9}\cap H_{16}\cap H_{20}\\
\nu_{12}=H_{1}\cap H_{6}\cap H_{10}\cap H_{14}\cap H_{18}.
\end{array}
\end{equation}
Let us now apply the generalized Delzant construction. 
It is easy to see that the following relations
$$
\left(
\begin{array}{c}
V_1\\
V_2\\
V_3
\end{array}
\right) =-\frac{1}{2}\left(\begin{array}{ccc}
\phi+\frac{1}{\phi}&\phi&\frac{1}{\phi}\\
\frac{1}{\phi}&\phi+\frac{1}{\phi}&\phi\\
\phi&\frac{1}{\phi}&\phi+\frac{1}{\phi}
\end{array}
\right) \left(\begin{array}{c} V_4\\V_5\\V_6\end{array}\right)
$$
$$
\left(
\begin{array}{c}
V_8\\
V_9\\
V_{10}
\end{array}
\right) =\frac{1}{2}\left(\begin{array}{ccc}
\frac{1}{\phi}&-1&-\phi\\
-\phi&\frac{1}{\phi}&-1\\
-1&-\phi&\frac{1}{\phi}
\end{array}
\right) \left(\begin{array}{c} V_4\\V_5\\V_6\end{array}\right)
$$
imply that the kernel of $\pi$, $\n$, is the $17$--dimensional subspace of $\R^{20}$ that is spanned by the vectors 
$$
\begin{array}{lll}
e_i&+&e_{10+i},\;\;i=1,\ldots,10\\
e_1&+&\frac{1}{2}[(\phi+\frac{1}{\phi})e_4+\phi e_5+\frac{1}{\phi}e_6]\\
e_2&+&\frac{1}{2}[\frac{1}{\phi}e_4+(\phi+\frac{1}{\phi})e_5+\phi e_6]\\
e_3&+&\frac{1}{2}[\phi e_4+\frac{1}{\phi}e_5+(\phi+\frac{1}{\phi}) e_6]\\
e_7&+&e_4+e_5+e_6\\
e_8&+&\frac{1}{2}[-\frac{1}{\phi}e_4+e_5+\phi e_6]\\
e_9&+&\frac{1}{2}[\phi e_4-\frac{1}{\phi}e_5+e_6]\\
e_{10}&+&\frac{1}{2}[e_4+\phi e_5-\frac{1}{\phi}e_6].\\
\end{array}
$$
Since the vectors $X_i$, $i=1,\ldots,20$, generate the quasilattice $B$, the group 
$N$ is connected and given by the group $\exp(\n)$.
The level set $\Psi^{-1}(0)$ for the moment mapping  $\Psi$ with respect to the induced $N$--action on $\C^{20}$ is given by the compact
subset of points $\underline{z}$ in $\C^{20}$ such that
$$
\left\{
\begin{array}{llll}
|z_i|^2&+&|z_{10+i}|^2=2&i=1,\ldots,10\\
|z_1|^2&+&\frac{1}{2}[(\phi+\frac{1}{\phi})|z_4|^2+\phi |z_5|^2+\frac{1}{\phi}|z_6|^2]=2\phi&\\
|z_2|^2&+&\frac{1}{2}[\frac{1}{\phi}|z_4|^2+(\phi+\frac{1}{\phi})|z_5|^2+\phi |z_6|^2]=2\phi&\\
|z_3|^2&+&\frac{1}{2}[\phi |z_4|^2+\frac{1}{\phi}|z_5|^2+(\phi+\frac{1}{\phi}) |z_6|^2]=2\phi&\\
|z_7|^2&+&|z_4|^2+|z_5|^2+|z_6|^2=4&\\
|z_8|^2&+&\frac{1}{2}[-\frac{1}{\phi}|z_4|^2+|z_5|^2+\phi |z_6|^2]=2&\\
|z_9|^2&+&\frac{1}{2}[\phi |z_4|^2 -\frac{1}{\phi}|z_5|^2+|z_6|^2]=2&\\
|z_{10}|^2&+&\frac{1}{2}[|z_4|^2+\phi |z_5|^2-\frac{1}{\phi}|z_6|^2]=2.&
\end{array}
\right.
$$
By Theorem~\ref{nonsemplice-re}, the quotient $M_{\Delta}=\Psi^{-1}(0)/N$ is stratified by symplectic quasifolds.

From the complex viewpoint, consider the open subset $\C^{20}_{\Delta}$ of $\C^{20}$ obtained by (\ref{aperto}) using
(\ref{tabella vertici icosaedro}). By Theorem~\ref{nonsemplice-cx}, the quotient 
$X_{\Delta}=\C_{\Delta}^{20}//N_{\SC}$ is stratified by complex quasifolds.

By Theorem~\ref{non semplice iso}, $M_{\Delta}$ and $X_{\Delta}$ can be identified. Their global description is similar to the one given for the symplectic and complex toric spaces corresponding to the octahedron, except that here $D^3=\R^3/B$ is a quasitorus and the regular stratum is a quasifold.

\begin{figure}[h]
\begin{center}
\includegraphics[height=60mm]{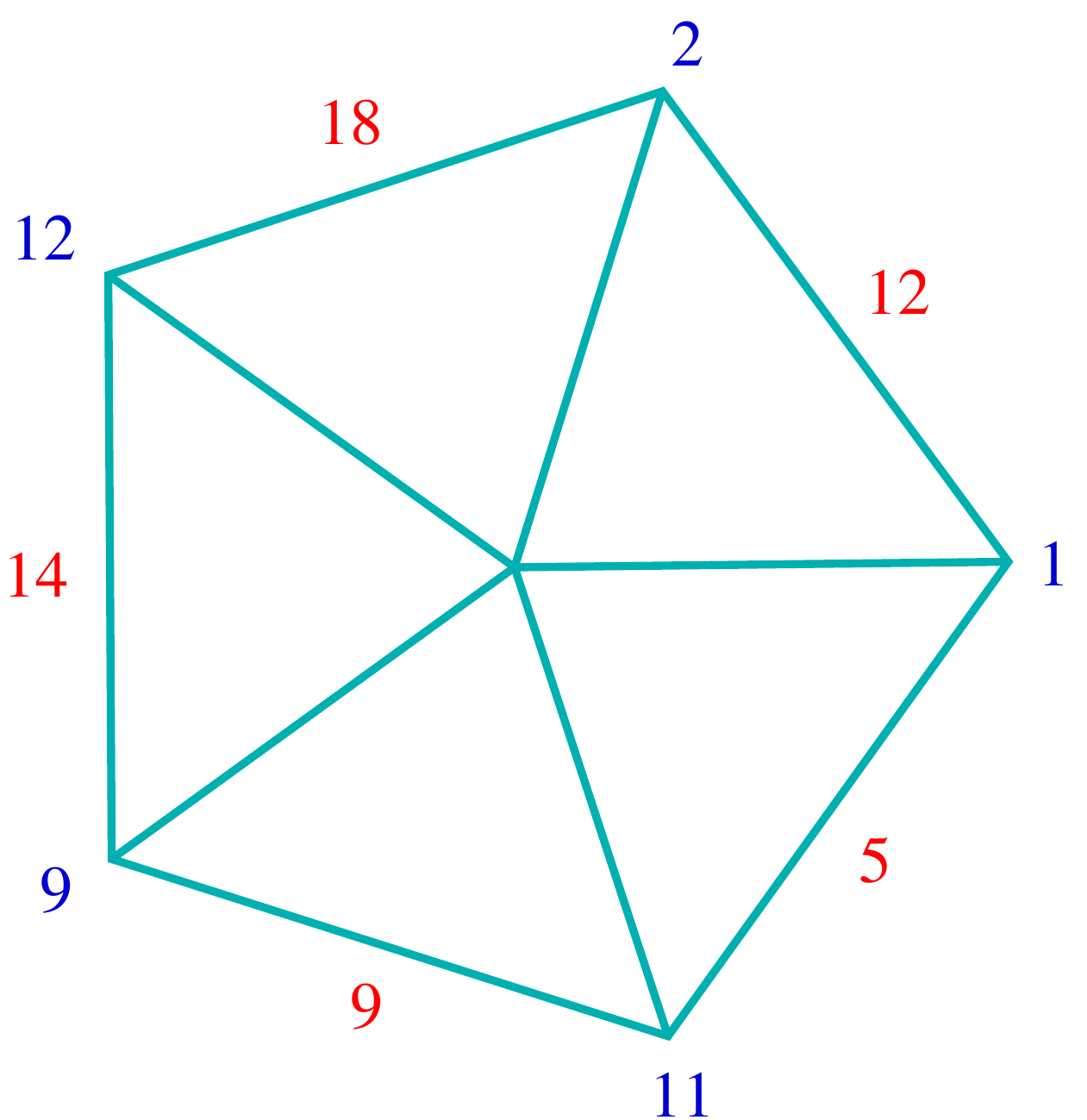}
\end{center}
\caption{The pentagon $\Delta_L$; its vertices and facets are suitably labeled.}\label{pentagono}
\end{figure}

We describe a chart for the regular part of the complex quotient. We choose the vertex $\nu_4$ and the triple $\{5,9,14\}$. A local model for the following open subset of the regular part $$\{\vz\in\C_{\Delta}^{20}\;|\;z_i\neq 0, i\neq 5,9,14\;\hbox{and}\;(z_5,z_9,z_{14})\neq(0,0,0)\}/N_{\SC}$$
is given by $(\C^3\smallsetminus\{0\})/\Gamma$,
where
$$\Gamma=\exp\{(-\phi h,\phi(h+k+2l),-\phi k)\;|\;h,k,l\in\Z\}.$$
Similarly to what happens for the octahedron, the quotient $X_{\Delta}$ around the singular point corresponding to $\nu_4$ is a cone.  
We will also have a diagram similar to (\ref{diagramma reale}). In this case 
$$\n({\cal C})=\{(-\phi s-t,\phi (s+t),s,-s-\phi t,t)\;|\;s,t\in\R\}$$
and
$$
N({\cal C})=\exp\{(-\phi s-t,\phi(s+t+2l),s,-s-\phi t,t)\;|\;s,t\in\R,\; l\in\Z\}.$$
The components of the moment mapping $\Psi\,\colon\C^5\longrightarrow(\n({\cal C}))^*$ with respect to the basis of
$(\n({\cal C}))^*$ dual to $(-\phi,\phi,1,-1,0),(-1,\phi,0,-\phi,1)$ are
$$
-\phi|z_1|^2+\phi|z_2|^2+|z_3|^2-|z_4|^2,\quad
-|z_1|^2+\phi|z_2|^2-|z_4|^2+|z_5|^2.
$$
As in the case of the octahedron, we can find the symplectic and complex toric spaces corresponding to the cone $\cal C$; these are regular symplectic and complex quasifolds away from the point corresponding to the cone apex.
We now describe the local structure of the stratification. Consider the vector 
$$X_{\nu_4}=2(-1,\phi,0)=-\frac{2}{2+\phi}(X_5+X_9+X_{12}+X_{14}+X_{18})=-V_1+V_2+V_3+V_4-V_5+V_6\in B.$$
The plane 
$H_{\nu_4}=\left\{\;\mu\in(\R^3)^*\;|\;\langle\mu,X_{\nu_4}\rangle=2/\phi\;\right\}$
cuts the icosahedron in the regular pentagon $\Delta_L$ whose vertices are given by $\{\nu_1,\nu_2,\nu_9,\nu_{11},\nu_{12}\}$  (see Figure~\ref{pentagono}). 
Let $\xi_{\nu_4}=(-\frac{1}{\phi(2+\phi)},\frac{1}{2+\phi},0)$ be a point in the cutting plane $H_{\nu_4}$.
The rigid motion $R+\xi_{\nu_4}\colon(\R^3)^*\rightarrow(\R^3)^*$, where $R$ is the rotation
$$\left(
\begin{array}{ccc}
\frac{\phi}{\sqrt{2+\phi}}&-\frac{1}{\sqrt{2+\phi}}&0\\
\frac{1}{\sqrt{2+\phi}}&\frac{\phi}{\sqrt{2+\phi}}&0\\
0&0&1
\end{array}
\right),
$$
maps the $xz$--plane into $H_{\nu_4}$. We obtain
$$\Delta_L=\bigcap_{i=5, 9,12,14,18}\{\;\mu=(\mu_1,0,\mu_3)\in(\R^3)^*\;|\;\langle R(\mu)+\xi_{\nu_4},X_{i}\rangle\geq-1\}.$$
This gives
$$\Delta_L=\bigcap_{i=5, 9,12,14,18}\left\{\;(\mu_1,\mu_3)\in(\R^2)^*\;|\;\langle(\mu_1,\mu_3),Y_i\rangle\geq-\frac{2}{2+\phi}\right\},$$
with 
$$\begin{array}{lcl}
Y_{5}&=&(-\frac{\phi}{\sqrt{2+\phi}},
\frac{1}{\phi})\\ 
Y_{9}&=&(\frac{1}{\phi\sqrt{2+\phi}},1)\\
Y_{12}&=&(-\frac{\phi}{\sqrt{2+\phi}},
 -\frac{1}{\phi})\\ Y_{14}&=&(\frac{2}{\sqrt{2+\phi}},0)\\ 
 Y_{18}&=&(\frac{1}{\phi\sqrt{2+\phi}},-1).
 \end{array}
$$ 
The above vectors are obtained from $X_5,X_9,X_{12},X_{14},X_{18}$ by rotating with $R^t$ and then projecting onto the
$xz$--plane. They generate the pentagonal quasilattice $B_L$ given by rotating and projecting the
icosahedral quasilattice $B$ onto the $xz$--plane. In fact, it is easy to check that $X_5,X_9,X_{12},X_{14},X_{18},X_{\nu_4}$ generate $B$ and,
moreover, that $R^t(X_{\nu_4})$ is parallel to $(0,1,0)$. 
Now we apply the generalized Delzant procedure to $\Delta_L$ with respect to the quasilattice $B_L$ and the vectors $Y_i$ above. 
Again here $N(\Delta_L)$ is connected and is given by
$\exp(\{(-\phi r-\phi s-t,r+\phi s+\phi t,s,-\phi r-s-\phi t,t)\;|\;r,s,t\in \R\})$. 
The level set $\Psi^{-1}(0)$ for the moment mapping  $\Psi$ with respect to the induced $N(\Delta_L)$--action on $\C^{5}$ 
is therefore the compact submanifold of $\C^5$ 
given by the points $\underline{z}$ in $\C^{5}$ such that
$$\left\{
\begin{array}{crcrcrclcl}
&\phi|z_1|^2&-&|z_2|^2&+&\phi|z_4|^2&&&=&\frac{2}{\phi}\\
-&\phi|z_1|^2&+&\phi|z_2|^2&+&|z_3|^2&-&|z_4|^2&=&0\\
-&|z_1|^2&+&\phi|z_2|^2&-&\phi|z_4|^2&+&|z_5|^2&=&0.
\end{array}
\right.
$$
We thus obtain the $6$--dimensional symplectic quasifold $M_{\Delta_L}=\Psi^{-1}(0)/N(\Delta_L)$. From the complex viewpoint, the toric quasifold $X_{\Delta_L}$ is given by the quotient $\C^5_{\Delta_L}/N_{\SC}(\Delta_L)$, where the open subset $\C^5_{\Delta_L}$ can be easily determined from Figure~\ref{pentagono}. 
We know that $N({\cal C})$ is a subgroup of
$N(\Delta_L)$; in fact 
$$N({\cal C})=\left\{\exp(-\phi r-\phi s-t,r+\phi s+\phi t,s,-\phi r-s-\phi t,t)\in N(\Delta_L)\;|\;r=-\frac{2}{\phi}l, \hbox{with}\,l\in\Z\right\}.$$
Interestingly, the quotient $N(\Delta_L)/N({\cal C})$ is the quasitorus $\R/2\phi\Z$. 
Thus, in this case, the fibers of the projections $p$ and $p_{\SC}$ of diagram (\ref{diagramma reale}) are given by $$N_{\SC}(\Delta_L)/N_{\SC}({\cal C})\simeq \R_{>0}\times \R/2\phi\Z.$$

\section*{Conflict of Interest}
The authors declare that there is no conflict of interests regarding the publication of this paper.

\section*{Acknowledgements}
This research was partially supported by  grant  PRIN 2010NNBZ78\_$\!$\_012 
(MIUR, Italy). All 3D pictures were drawn using ZomeCAD.

\end{document}